\documentclass[12pt]{amsart} 
\usepackage{amssymb,amsmath,amscd} 

\def\mysavedown#1{\edef\mysubs{\mysubs#1}}
\def\mysaveup#1{\edef\mysups{\mysups#1}}
\def\mydown#1{{\mytensor}_{\vphantom{\mysubs}#1}}
\def\myup#1{{\mytensor}^{\vphantom{\mysups}#1}}
\def\tensor#1#2{
  #1
  \def\mytensor{\vphantom{#1}}
  \def\mysubs{\relax}
  \def\mysups{\relax}
  \let\down=\mysavedown
  \let\up=\mysaveup
  #2
  \let\down=\mydown
  \let\up=\myup
  #2
  }

\theoremstyle{plain} 
\newtheorem{Lem}{Lemma}[section] 
\newtheorem{Prop}[Lem]{Proposition} 
\newtheorem{Thm}[Lem]{Theorem} 
\newtheorem{Cor}[Lem]{Corollary} 
\theoremstyle{definition} 
\newtheorem{Def}[Lem]{Definition} 
\newtheorem{Rem}[Lem]{Remark} 
 
\newtheorem{Prbl}[Lem]{Problem} 
\newtheorem{Proc}[Lem]{Procedure} 
\errorcontextlines=0  \numberwithin{equation}{section}

\newcommand{\bbC}{{\mathbb C}} 
\newcommand{\bbR}{{\mathbb R}}

\newcommand{\bbN}{{\mathbb N}} 
\newcommand{\bbK}{{\mathbb K}} 
\newcommand{\calA}{{\mathcal A}} 
\newcommand{\calB}{{\mathcal B}} 
\newcommand{\calH}{{\mathcal H}} 
\newcommand{\calL}{{\mathcal L}}

\newcommand{\calP}{{\mathcal P}} 
\newcommand{\calS}{{\mathcal S}} 
\newcommand{\calT}{{\mathcal T}} 
\newcommand{\calV}{{\mathcal V}} 
 
\newcommand{\id}{\mathrm{id}} 
\newcommand{\frakl}{{\mathfrak l}} 
\newcommand{\frakr}{{\mathfrak r}} 

\newcommand{\fR}{{\mathfrak R}}
 
\newcommand{\bpi}{\begin{picture}} 
\newcommand{\epi}{\end{picture}}

\newcommand{\twoline}[2]{\genfrac{}{}{0pt}{}{#1}{#2}}

\begin{document} 

\title[Algebraic covariant derivative curvature tensors] 
{Short formulas for algebraic covariant derivative curvature tensors via Algebraic Combinatorics} 
\author[B. Fiedler]{Bernd Fiedler}
\address{Bernd Fiedler \\ Mathematisches Institut \\ Universit\"at Leipzig\\ 
Augustusplatz 10/11 \\ D-04109 Leipzig \\ Germany}
\urladdr{http://home.t-online.de/home/Bernd.Fiedler.RoschStr.Leipzig/}  
\email{Bernd.Fiedler.RoschStr.Leipzig@t-online.de}  
\subjclass{53B20, 15A72, 05E10, 16D60, 05-04} 

\begin{abstract}
We consider generators of algebraic covariant derivative curvature tensors ${\mathfrak R}'$ which can be constructed by a Young symmetrization of product tensors $W\otimes U$ or $U\otimes W$, where $W$ and $U$ are covariant tensors of order 2 and 3. $W$ is a symmetric or alternating tensor whereas $U$ belongs to a class of the infinite set ${\mathfrak S}$ of irreducible symmetry classes characterized by the partition $(2\,1)$. Using Computer Algebra we search for such generators whose coordinate representations are polynomials with a minimal number of summands. For a generic choice of the symmetry class of $U$ we obtain lengths of 16 or 20 summands if $W$ is symmetric or skew-symmetric, respectively. In special cases these numbers can be reduced to the minima 12 or 10. If these minima occur then $U$ admits an index commutation symmetry. Furthermore minimal lengths are possible if $U$ is formed from torsion-free covariant derivatives of symmetric or alternating 2-tensor fields.

Foundation of our investigations is a theorem of S. A. Fulling, R. C. King, B. G. Wybourne and C. J. Cummins about a Young symmetrizer that generates the symmetry class of algebraic covariant derivative curvature tensors.
Furthermore we apply ideals and idempotents in group rings ${\mathbb C}[{\calS}_r]$ and discrete Fourier transforms for symmetric groups ${\calS}_r$. For symbolic calculations we used the Mathematica packages {\sf Ricci} and {\sf PERMS}.
\end{abstract}

\maketitle 

%
%

\section{Introduction}

The present paper continues investigations of \cite{fie20,fie03b} in which we constructed generators of {\itshape algebraic curvature tensors} and {\itshape algebraic covariant derivative curvature tensors}.

Algebraic curvature tensors are tensors of order 4 which have the same symmetry properties as the Riemann tensor of a Levi-Civita connection in Differential Geometry.
Let ${\calT}_r V$ be the vector space of the $r$-times covariant tensors $T$ over a finite-dimensional $\bbK$-vector space $V$, $\bbK = \bbR$ or  $\bbK = \bbC$. We assume that $V$ possesses a {\itshape fundamental tensor} $g \in {\calT}_2 V$ (of arbitrary signature) which can be used for raising and lowering of tensor indices.
\begin{Def}
A tensor $\fR \in {\calT}_4 V$ is called an {\itshape algebraic curvature tensor} iff $\fR$ has the {\itshape index commutation symmetry}
\begin{eqnarray}
\forall\,w, x, y, z\in V:\;\;\;
\fR(w,x,y,z) & = & - \fR(w,x,z,y) \;=\; \fR(y,z,w,x)
\label{equ1.1}
\end{eqnarray}
and fulfills the {\itshape first Bianchi identity}
\begin{eqnarray}
\forall\,w, x, y, z\in V:\;\;\;\fR(w,x,y,z) + \fR(w,y,z,x) + \fR(w,z,x,y) & = & 0\,.
\end{eqnarray}
\end{Def}
\begin{Def}
A tensor $\fR' \in {\calT}_5 V$ is called an {\itshape algebraic covariant derivative curvature tensor} iff $\fR'$ has the {\itshape index commutation symmetry}
\begin{eqnarray}
\fR'(w,x,y,z,u) & = & - \fR'(w,x,z,y,u) \;=\; \fR'(y,z,w,x,u)
\end{eqnarray}
and fulfills the {\itshape first Bianchi identity}
\begin{eqnarray}
\fR'(w,x,y,z,u) + \fR'(w,y,z,x,u) + \fR'(w,z,x,y,u) & = & 0
\end{eqnarray}
and the {\itshape second Bianchi identity}
\begin{eqnarray}
\fR'(w,x,y,z,u) + \fR'(w,x,z,u,y) + \fR'(w,x,u,y,z) & = & 0 \label{equ1.5}
\end{eqnarray}
for all $u, w, x, y, z\in V$.
\end{Def}
The relations (\ref{equ1.1}) -- (\ref{equ1.5}) correspond to the well-known formulas
\begin{eqnarray}
 & & R_{i j k l}\;=\; - R_{i j l k}\;=\;R_{k l i j} \label{equ1.6}\\
 & & R_{i j k l} + R_{i k l j} + R_{i l j k}\;=\;0
\end{eqnarray}
for the Riemann tensor $R$ and
\begin{eqnarray}
 & & R_{i j k l\,;\,m}\;=\; - R_{i j l k\,;\,m}\;=\;R_{k l i j\,;\,m}\\
 & & R_{i j k l\,;\,m} + R_{i k l j\,;\,m} + R_{i l j k\,;\,m}\;=\;0\\
 & & R_{i j k l\,;\,m} + R_{i j l m\,;\,k} + R_{i j m k\,;\,l}\;=\;0
\end{eqnarray}
for its first covariant derivative which we present here in terms of tensor coordinates with respect to an arbitrary local coordinate system.

A famous problem connected with algebraic curvature tensors is the {\itshape Osserman conjecture}.
\begin{Def}
Let $\fR \in {\calT}_4 V$ be an algebraic curvature tensor. For $x \in V$, the {\itshape Jacobi operator} $J_{\fR}(x)$ of $\fR$ and $x$ is the linear operator
$J_{\fR}(x) : V \rightarrow V\;,\; J_{\fR}(x): y \mapsto J_{\fR}(x)y$ that is defined by
$\forall\, w \in V :\; g(J_{\fR}(x) y , w) = \fR(y, x, x, w)$.
\end{Def}
%
%
%
\begin{Def}
An algebraic curvature tensor $\fR$ is called {\itshape Osserman} if the eigenvalues of $J_{\fR}(x)$ are constant both on $S^{+}(V):=\{x\in V\,|\,g(x,x)=+1\}$ and on $S^{-}(V):=\{x\in V\,|\,g(x,x)=-1\}$.
\end{Def}
If $R$ is the Riemann tensor of a Riemannian manifold $(M,g)$ which is locally a rank one symmetric space or flat, then the eigenvalues of $J_R(x)$ are constant on the unit sphere bundle of $(M,g)$. Osserman \cite{oss90} wondered if the converse held. This question is known as the {\itshape Osserman conjecture}.

The correctness of the Osserman conjecture has been established for Riemannian manifolds $(M,g)$ in all dimensions $\not= 8, 16$ (see \cite{c88,n02}) and for Lorentzian manifolds $(M,g)$ in all dimensions (see \cite{bbg97,gkv97}). However Osserman's question has a negative answer in the case of a pseudo-Riemannian metric with signature $(p,q)$, $p,q\ge 2$ (see e.g. the references given in \cite[p.2]{fie20}).
A detailed view about the Osserman conjecture can be found in the book \cite{gilkey5} by P. B. Gilkey.

Numerous examples of Osserman algebraic curvature tensors can be constructed by menas of operators $\alpha$ and $\gamma$ given below. It turned out that these operators lead to generators for arbitrary algebraic curvature tensors.
\begin{Def}
Let $S , A \in {\calT}_2 V$ be symmetric or alternating tensors of order 2, i.e. their coordinates satisfy $S_{i j} = S_{j i}$, $A_{i j} = - A_{j i}$. We define tensors $\gamma (S), \alpha(A) \in {\calT}_4 V$ by
\begin{eqnarray}
\gamma (S)_{i j k l} & := &
{\textstyle \frac{1}{3}} \left( S_{i l} S_{j k} - S_{i k} S_{j l} \right) \,,
\label{equ1.12} \\
\alpha (A)_{i j k l} & := &
{\textstyle \frac{1}{3}} \left( 2\,A_{i j} A_{k l} + A_{i k} A_{j l}
-  A_{i l} A_{j k} \right) \,. \label{equ1.13}
\end{eqnarray}
\end{Def}
Now we can construct an example of an Osserman algebraic curvature tensor in the following way. Let $g\in{\calT}_2 V$ be a positive definite metric and $\{ C_i \}_{i = 1}^r$ be a finite set of real, skew-symmetric
$(\mathrm{dim} V \times \mathrm{dim} V)$-matrices that satisfy the Clifford commutation relations
$C_i \cdot C_j + C_j \cdot C_i = - 2\,{\delta}_{i j}\,{\rm Id}$.
If we form skew-symmetric tensors $A_i\in\calT_2 V$ by $A_i{(x,y)} := g(C_i\cdot x,y)$ ($x,y\in V$), then
\begin{eqnarray}
\fR & = & {\lambda}_0\,\gamma (g) + \sum_{i = 1}^r
{\lambda}_i \,\alpha (A_i) \;\;\;,\;\;\;{\lambda}_0 , {\lambda}_i = \mathrm{const.}
\end{eqnarray}
is an Osserman algebraic curvature tensor (see \cite{gilkey3}). Further examples which allow also indefinite metrics can be found in \cite[pp.191-193]{gilkey5}. (See also \cite[Sec.6]{fie20}.)

The operators $\alpha$ and $\gamma$ can be used to form generators for arbitrary algebraic curvature tensors. P. Gilkey \cite[pp.41-44]{gilkey5} and B. Fiedler \cite{fie20} gave different proofs for
\begin{Thm}
Each of the sets of tensors
\begin{enumerate}
\item{$\{ \gamma(S)\;|\; S\in\calT_2 V\;\mathrm{symmetric} \}$}
\item{$\{ \alpha(A)\;|\; A\in\calT_2 V\;\mbox{\rm skew-symmetric} \}$}
\end{enumerate}
generates the vector space of all algebraic curvature tensors $\fR$ on $V$.
\end{Thm}
Note that the tensors $\gamma(S)$ and $\alpha(A)$ are expressions which arise from
$S\otimes S$ or $A\otimes A$ by a symmetrization
$\gamma(S) = \textstyle{\frac{1}{12}}\,y_t^{\ast}(S\otimes S)$,
$\alpha(A) = \textstyle{\frac{1}{12}}\,y_t^{\ast}(A\otimes A)$,
where $y_t$ is the Young symmetrizer of the Young tableau
\begin{eqnarray}
t & = 
\begin{array}{|c|c|}
\hline
1 & 3 \\
\hline
2 & 4 \\
\hline
\end{array}\,.
\end{eqnarray}
(See \cite{fie20}. See also Section 2 for basic facts and definitions.)

In \cite{fie03b} we searched for similar {\itshape generators of algebraic covariant derivative curvature tensors}. We used Boerner's definition of {\itshape symmetry classes} for tensors $T\in\calT_r V$ by right ideals $\frakr\subseteq\bbK[\calS_r]$ of the group ring $\bbK[\calS_r]$ of the symmetric group $\calS_r$ (see Section \ref{sec2} and \cite{boerner,boerner2,fie16,fie18}). On this basis we investigated the following
\begin{Prbl}
We search for generators of algebraic covariant derivative curvature tensors which can be formed by a suitable symmetry operator from tensors
\begin{eqnarray}
T\otimes \hat{T}\;\;\;\mathrm{or}\;\;\;\hat{T}\otimes T &\;\;\;,\;\;\;&
T\in\calT_2 V\;,\;\hat{T}\in\calT_3 V \label{equ1.18}
\end{eqnarray}
where $T$ and $\hat{T}$ belongs to symmetry classes of $\calT_2 V$ and $\calT_3 V$ which are defined by minimal right ideals $\frakr\subset\bbK[\calS_2]$ and $\hat{\frakr}\subset\bbK[\calS_3]$, respectively. \label{probl1.7}
\end{Prbl}
All such generators can be gained by means of the {\itshape Young symmetrizer} $y_{t'}$ of the {\itshape Young tableau}
\begin{eqnarray}
t' & = &
\begin{array}{|c|c|c|c}
\cline{1-3}
1 & 3 & 5 & \\
\cline{1-3}
2 & 4 &\multicolumn{2}{c}{\;\;\;} \\
\cline{1-2}
\end{array}
\,.
\end{eqnarray}
In \cite{fie03b} we proved the following Theorems \ref{thm1.8}--\ref{thm1.10}:
\begin{Thm} \label{thm1.8}
A solution of Problem {\rm\ref{probl1.7}} can be constructed at most from such pairs of tensors {\rm (\ref{equ1.18})} whose symmetry classes are characterized by the following partitions $\lambda\vdash 2$, $\hat{\lambda}\vdash 3$:
\begin{center}
{\rm
\begin{tabular}{|c|cl|cl|}
\hline
 & \multicolumn{2}{c|}{$\lambda$ for $T$} & \multicolumn{2}{c|}{$\hat{\lambda}$ for $\hat{T}$} \\
\hline
(a) & $(2)$, & i.e. $T$ symmetric & $(3)$, & i.e. $\hat{T}$ symmetric \\
(b) & $(2)$, & i.e. $T$ symmetric & $(2\,1)$ & \\
(c) & $(1^2)$, & i.e. $T$ skew-symmetric & $(2\,1)$ & \\
\hline
\end{tabular}\,.
}
\end{center}
\end{Thm}
The proof of Theorem \ref{thm1.8} is based on the {\it Littlewood-Richardson rule}\footnote{See the
references \cite{kerber,kerber3,jameskerb,littlew1,mcdonald,full4,fultharr} for the Littlewood-Richardson rule.}.
 (see \cite{fie03b}). The case (a) of Theorem \ref{thm1.8} is specified by
\begin{Thm}
Let us denote by $S\in \calT_2 V$ and $\hat{S}\in \calT_3 V$ symmetric tensors of order {\rm 2} and {\rm 3}, respectively.
Then the set of all tensors which belong to exactly one of the following tensor types
\begin{eqnarray}
\tau: & &
\begin{array}{cccc}
y_{t'}^{\ast} (S\otimes \hat{S}) & , & y_{t'}^{\ast} (\hat{S}\otimes S) & , \\
\end{array} \label{equ1.20}
\end{eqnarray}
generates the vector space of all algebraic covariant derivative curvature tensors $\fR'\in\calT_5 V$.
Moreover, the tensors {\rm (\ref{equ1.20})} coincide and their coordinates fulfill
\begin{eqnarray} \label{equ1.21}
\hat{\gamma}(S,\hat{S})_{ijkls} & := &
(y_{t'}^{\ast} (S\otimes \hat{S}))_{ijkls} \;=\; (y_{t'}^{\ast} (\hat{S}\otimes S))_{ijkls} \\
 & = & 4\,\left\{S_{il}{\hat{S}}_{jks} - S_{jl}{\hat{S}}_{iks} + S_{jk}{\hat{S}}_{ils} - S_{ik}{\hat{S}}_{jls} \right\} \nonumber
\end{eqnarray} \label{thm1.9}
\end{Thm}
The operator $\hat{\gamma}$ plays the same role for the generators of algebraic covariant derivative curvature tensors considered in Theorem \ref{thm1.9} as the operators $\alpha$ and $\gamma$ play for the generators of algebraic curvature tensors. A first proof that the expressions (\ref{equ1.21}) are generators for $\fR'$ was given by P. B. Gilkey \cite[p.236]{gilkey5}.

The cases (b) and (c) of Theorem \ref{thm1.8} lead to
\begin{Thm} \label{thm1.10}
Let us denote by $S, A \in \calT_2 V$ symmetric or alternating tensors of order {\rm 2} and by $U \in \calT_3 V$ covariant tensors of order {\rm 3} whose symmetry class $\calT_{\frakr}$ is defined by a fixed minimal right ideal $\frakr\subset\bbK[\calS_3]$ from the equivalence class characterized by the partition $(2 , 1) \vdash 3$. We consider the following types $\tau$ of tensors
\begin{eqnarray} \label{equ1.22}%
\tau: & &
y_{t'}^{\ast} (S \otimes U)\;\;,\;\;y_{t'}^{\ast} (U \otimes S)\;\;,\;\;
y_{t'}^{\ast} (A \otimes U)\;\;,\;\;y_{t'}^{\ast} (U \otimes A)\,.
\end{eqnarray}
Then for each of the  above types $\tau$ the following assertions are equivalent:
\begin{enumerate}
\item{The vector space of algebraic covariant derivative curvature tensors $\fR' \in \calT_5 V$ is the set of all finite sums of tensors of the type $\tau$ considered. \label{statement1}}
\item{The right ideal $\frakr$ is different from the right ideal
$\frakr_0 := f \cdot \bbK [\calS_3]$ with generating idempotent
\begin{eqnarray}
f \;:=\; \left\{ \frac{1}{2}\,(\id - (1 \,3)) - \frac{1}{6}\,y \right\}
& \;\;\;,\;\;\; &
y \;:=\; \sum_{p \in \calS_3} \mathrm{sign}(p)\,p \,. \label{equ1.23}
\end{eqnarray}
}
\end{enumerate}
\end{Thm}
In Theorem \ref{thm1.10} not only a single symmetry class is allowed for the tensors $U$ but the {\it complete infinite set} ${\mathfrak S}$ of irreducible symmetry classes of $(2\,1)\vdash 3$ from which only the class of the right ideal $\frakr_0$ has to be excluded.

In the situation of Theorem \ref{thm1.10} we can also determine operators of the type $\alpha$, $\gamma$, $\hat{\gamma}$ which describe the generators of the algebraic covariant derivative curvature tensors $\fR'$ considered. However, these operators depend on the right ideal $\frakr$ (or its generating idempotent $e$) that defines the symmetry class of $U$. And they yield no short expressions of 2, 3, or 4 terms but longer expressions between 10 and 20 terms of length. {\it The search for shortest expressions of this type is the subject of our paper.} Some of our main results are collected in
\begin{Thm} \label{thm1.11}%
Consider the situation of Theorem {\rm\ref{thm1.10}}. Then it holds:
\begin{enumerate}
\item{The tensors {\rm (\ref{equ1.22})} satisfy
\begin{eqnarray}
y_{t'}^{\ast} (S \otimes U)\;=\;y_{t'}^{\ast} (U \otimes S) & \;,\;
y_{t'}^{\ast} (A \otimes U)\;=\; - y_{t'}^{\ast} (U \otimes A)\,.
\end{eqnarray}
}
\item{Let $\dim V\ge 3$. Then the coordinates of {\rm (\ref{equ1.22})} are sums of the following lengths
\begin{center}
{\rm
\begin{tabular}{|c|l|c|c|}
\hline
 & & $y_{t'}^{\ast} (S \otimes U)$ & $y_{t'}^{\ast} (A \otimes U)$ \\
\hline
(a) & generic case for $\frakr$ & 16 & 20 \\
(b) & $\frakr$ producing minimal lengths & 12 & 10 \\
\hline
\end{tabular}\,.
}
\end{center}
}
\item{For $\dim V\ge 3$ there exist minimal right ideals $\frakr$ of $(2\,1)\vdash 3$ which lead to the minimal lengths of case {\rm (b)} both for $y_{t'}^{\ast} (S \otimes U)$ and for $y_{t'}^{\ast} (A \otimes U)$.}
\item{If the coordinates of $y_{t'}^{\ast} (S \otimes U)$ or $y_{t'}^{\ast} (A \otimes U)$ have the minimal lengths of case {\rm (b)} and $\dim V\ge 3$ then $U$ admits an index commutation symmetry.}
\end{enumerate}
\end{Thm}
Further results are given in Section \ref{sec5} and the Appendices of our paper.

The concept ''expression of minimal length'' depends on the method which we use to reduce expressions (see Section \ref{sec4}). Remark \ref{rem4.9} discusses a generalization of our reduction method which could possibly lead to a further decrease of the numbers in Theorem \ref{thm1.11}.

Examples of tensors $U$ with a symmetry from ${\mathfrak S}$ are given by tensor fields
\begin{eqnarray} \label{equ1.22a}%
U\;=\;(\nabla\psi - {\rm sym}(\nabla\psi))|_p & {\rm or} &
U\;=\;(\nabla\omega - {\rm alt}(\nabla\omega))|_p\;\;\;,\;\;\;p\in M\,,
\end{eqnarray}
where $\psi, \omega\in\calT_2 M$ are symmetric/alternating tensor fields of order 2 on a diffentiable manifold $M$ and $\nabla$ is a torsion-free covariant derivative (see \cite{fie03a} and Section \ref{sec6}). For tensors (\ref{equ1.22a}) we obtained the following result
\begin{Thm} \label{thm1.12}%
If we consider tensors $U$, $S$, $A$ on a tangent space $V = T_p M$ of a differentiable manifold $M$, $\dim M\ge 3$, and generate $U$ by one of the formulas {\rm (\ref{equ1.22a})} then we obtain the shortest lengths from Theorem {\rm \ref{thm1.11}, (2b)} exactly in the following cases:
\begin{enumerate}
\item{$y_{t'}^{\ast}(U\otimes S)$ and $U = (\nabla\psi - {\rm sym}(\nabla\psi))|_p$, $\psi\in\calT_2 M$ symmetric,}
\item{$y_{t'}^{\ast}(U\otimes S)$, $y_{t'}^{\ast}(U\otimes A)$ and
$U = (\nabla\omega - {\rm alt}(\nabla\omega))|_p$, $\omega\in\calT_2 M$ skew-symmetric.}
\end{enumerate}
\end{Thm}

Here is a brief outline to the paper. In Section \ref{sec2} we give a summary of basic facts about symmetry classes, Young symmetrizers and discrete Fourier transforms. These tools are needed to obtain the infinite set ${\mathfrak S}$ of symmetry classes for $U$. In Section \ref{sec4} and \ref{sec5} we construct short coordinate representations for the tensors (\ref{equ1.22}) by determining and solving a complete system of linear identities for the tensors $U$. In Section \ref{sec6} we determine the conditions for the occurence of index commutation symmetries on the tensors $U$ and show statement (4) of Theorem \ref{thm1.11}. In Section \ref{sec7} we prove Theorem \ref{thm1.12}.
Many results were obtained by computer calculations by means of the {\sf Mathematica} packages {\sf Ricci} \cite{ricci3} and {\sf PERMS} \cite{fie10}. The {\sf Mathematica} notebooks of these calculations are available at \cite{fie21}.\vspace{20pt}

\section{Basic facts} \label{sec2}
The vector spaces of algebraic curvature tensors or algebraic covariant derivative tensors over $V$ are {\itshape symmetry classes} in the sence of H. Boerner \cite[p.127]{boerner}. We denote by $\bbK [{\calS}_r]$ the {\itshape group ring} of a symmetric group ${\calS}_r$ over the field $\bbK$. Every group ring element $a = \sum_{p \in {\calS}_r} a(p)\,p \in \bbK [{\calS}_r]$ acts as so-called {\itshape symmetry operator} on tensors $T \in {\calT}_r V$ according to the definition
\begin{eqnarray}
(a T)(v_1 , \ldots , v_r) & := & \sum_{p \in {\calS}_r} a(p)\,
T(v_{p(1)}, \ldots , v_{p(r)}) \;\;\;\;\;,\;\;\;\;\;
v_i \in V \,. \label{equ2.1}
\end{eqnarray}
Equation \eqref{equ2.1} is equivalent to
$(a T)_{i_1 \ldots i_r} = \sum_{p \in {\calS}_r} a(p)\,
T_{i_{p(1)} \ldots  i_{p(r)}}$.
\begin{Def}
Let $\frakr \subseteq \bbK [{\calS}_r]$ be a right ideal of $\bbK [{\calS}_r]$ for which an $a \in \frakr$ and a $T \in {\calT}_r V$ exist such that $aT \not= 0$. Then the tensor set
\begin{eqnarray}
{\calT}_{\frakr} & := & \{ a T \;|\; a \in \frakr \;,\;
T \in {\calT}_r V \}
\end{eqnarray}
is called the {\itshape symmetry class} of tensors defined by $\frakr$.
\end{Def}
Since $\bbK [{\calS}_r]$ is semisimple for $\bbK = \bbR , \bbC$, every right ideal $\frakr \subseteq \bbK [{\calS}_r]$ possesses a generating idempotent $e$, i.e. $\frakr$ fulfils $\frakr = e \cdot \bbK [{\calS}_r]$. It holds (see e.g. \cite{fie20} or \cite{boerner,boerner2})
\begin{Lem}
If $e$ is a generating idempotent of $\frakr$, then a tensor $T \in {\calT}_r V$ belongs to ${\calT}_{\frakr}$ iff
$e T = T$.
Thus we have
${\calT}_{\frakr} = \{ eT \;|\; T \in {\calT}_r V \}$.
\end{Lem}
Now we summarize tools from our Habilitationsschrift \cite{fie17} (see also its summary \cite{fie18}).
We make use of the following connection
between $r$-times covariant tensors $T \in {\calT}_r V$
and elements of the {\it group ring}
${\bbK} [{\calS}_r]$.
\begin{Def} \label{def2.3}
 Any tensor
 $T \in {\calT}_r V$
 and any $r$-tuple
 $b := (v_1 , \ldots , v_r ) \in V^r$
 of
 $r$
 vectors from
 $V$
 induce a function
 $T_b : {\calS}_r \rightarrow {\bbK}$
 according to the rule
 \begin{eqnarray}
T_b (p) & := & T(v_{p(1)} , \ldots , v_{p(r)})\;\;\;,\;\;\;p \in {\calS}_r \,.
\end{eqnarray}
We identify this function with the group ring element
$T_b := \sum_{p \in {\calS}_r}T_b (p)\,p \in {\bbK} [{\calS}_r]$.
\end{Def}
Obviously,
two tensors $S , T \in {\calT}_r V$ fulfil $S = T$ iff
$S_b = T_b$ for all
$b \in V^r$.
We denote by '$\ast$'
the mapping
$\ast : a = \sum_{p \in {\calS}_r} a(p)\,p \;\mapsto\; a^{\ast} :=
\sum_{p \in {\calS}_r} a(p)\,p^{-1}$. Then the following important formula\footnote{See B. Fiedler \cite[Sec.III.1]{fie16} and B. Fiedler \cite{fie17}.} holds
\begin{eqnarray} \label{equ2.7}
\forall\,T\in\calT_r V\;,\;a\in\bbK[\calS_r]\;,\;b\in V^r\;:\;\;\;\;
(a\,T)_b & = & T_b\cdot a^{\ast}\,.
\end{eqnarray}
Now it can be shown that all $T_b$ of tensors $T$ of a given symmetry class lie in a certain left ideal of ${\bbK}[{\calS}_r]$.
\begin{Prop}\hspace{-1mm}\footnote{See B. Fiedler \cite{fie17} or
B. Fiedler \cite[Prop. III.2.5, III.3.1, III.3.4]{fie16}.}
\label{prop2.4}%
Let $e \in {\bbK}[{\calS}_r]$ be an idempotent. Then a
$T \in {\calT}_r V$ 
fulfils the condition
$eT = T$
iff
$T_b \;\in\; {\frakl} := {\bbK} [{\calS}_r] \cdot e^{\ast}$ for all
$b \in V^r$, i.e.
all $T_b$ of $T$
lie in the left ideal ${\frakl}$ generated by $e^{\ast}$.
\end{Prop}
The proof follows easily from (\ref{equ2.7}). Since a rigth ideal $\frakr$ defining a symmetry class and the left ideal $\frakl$ from Proposition \ref{prop2.4} satisfy $\frakr = \frakl^{\ast}$, we denote symmetry classes also by $\calT_{\frakl^{\ast}}$. A further result is
\begin{Prop}\hspace{-1mm}\footnote{See B. Fiedler \cite{fie17} or
B. Fiedler \cite[Prop. III.2.6]{fie16}.}
\label{prop2.5}%
If $\dim V \ge r$, then every left ideal
${\frakl} \subseteq {\bbK}[{\calS}_r]$ fulfils
${\frakl} = {\calL}_{\bbK} \{ T_b \;|\;
T \in {\calT}_{{\frakl}^{\ast}} \,,\, b \in V^r \}$.
(Here ${\calL}_{\bbK}$ denotes the forming of the linear closure.)
\end{Prop}
If $\dim V < r$, then the $T_b$ of the tensors from
${\calT}_{{\frakl}^{\ast}}$ 
will span only a linear subspace of
${\frakl}$ 
in general.

Important special symmetry operators are Young symmetrizers, which are defined by means of Young tableaux.

A {\itshape Young tableau} $t$ of $r\in\bbN$ is an arrangement of $r$ boxes such that
\begin{enumerate}
\item{the numbers ${\lambda}_i$ of boxes in the rows $i = 1 , \ldots , l$ form a decreasing sequence
${\lambda}_1 \ge {\lambda}_2 \ge \ldots \ge {\lambda}_l > 0$ with
${\lambda}_1 + \ldots + {\lambda}_l = r$,}
\item{the boxes are fulfilled by the numbers $1, 2, \ldots , r$ in any order.}
\end{enumerate}
For instance, the following graphics shows a Young tableau of $r = 16$.
\[\left.
\begin{array}{cc|c|c|c|c|c|c}
\cline{3-7}
{\lambda}_1 = 5 & \;\;\; & 11 & 2 & 5 & 4 & 12 & \\
\cline{3-7}
{\lambda}_2 = 4 & \;\;\; & 9 & 6 & 16 & 15 & \multicolumn{2}{c}{\;\;\;} \\
\cline{3-6}
{\lambda}_3 = 4 & \;\;\; & 8 & 14 & 1 & 7 & \multicolumn{2}{c}{\;\;\;} \\
\cline{3-6}
{\lambda}_4 = 2 & \;\;\; & 13 & 3 & \multicolumn{4}{c}{\hspace{2cm}} \\
\cline{3-4}
{\lambda}_5 = 1 & \;\;\; & 10 & \multicolumn{4}{c}{\hspace{2cm}} \\
\cline{3-3}
\end{array}\right\}\;=\;t\,.
\]
Obviously, the unfilled arrangement of boxes, the {\itshape Young frame}, is characterized by a partition
$\lambda = ({\lambda}_1 , \ldots , {\lambda}_l) \vdash r$ of $r$.

If a Young tableau $t$ of a partition $\lambda \vdash r$ is given, then the {\itshape Young symmetrizer} $y_t$ of $t$ is defined by\footnote{We use the convention $(p \circ q) (i) := p(q(i))$ for the product of two permutations $p, q$.}
\begin{eqnarray}
y_t & := & \sum_{p \in {\calH}_t} \sum_{q \in {\calV}_t} \mathrm{sign}(q)\, p \circ q
\end{eqnarray}
where ${\calH}_t$, ${\calV}_t$ are the groups of the {\itshape horizontal} or
{\itshape vertical permutations} of $t$ which only permute numbers within rows or columns of $t$, respectively. The Young symmetrizers of $\bbK [{\calS}_r]$ are essentially idempotent and define decompositions
\begin{eqnarray}
\bbK [{\calS}_r] \;=\;
\bigoplus_{\lambda \vdash r} \bigoplus_{t \in {\calS\calT}_{\lambda}}
\bbK [{\calS}_r]\cdot y_t
& \;\;,\;\; &
\bbK [{\calS}_r] \;=\;
\bigoplus_{\lambda \vdash r} \bigoplus_{t \in {\calS\calT}_{\lambda}}
y_t \cdot \bbK [{\calS}_r] \label{decomp}
\end{eqnarray}
of $\bbK [{\calS}_r]$ into minimal left or right ideals. In \eqref{decomp}, the symbol ${\calS\calT}_{\lambda}$ denotes the set of all standard tableaux of the partition $\lambda$. Standard tableaux are Young tableaux in which the entries of every row and every column form an increasing number sequence.\footnote{About Young symmetrizers and
Young tableaux see for instance
\cite{boerner,boerner2,full4,fulton,jameskerb,kerber,littlew1,mcdonald,%
muell,naimark,%
waerden,weyl1}. In particular, properties of Young symmetrizers in the case
${\bbK} \not= {\bbC}$ are described in \cite{muell}.}

S.A. Fulling, R.C. King, B.G.Wybourne and C.J. Cummins showed in \cite{full4} that the symmetry classes of the Riemannian curvature tensor $R$ and its {\itshape symmetrized\footnote{$(\,\ldots\,)$ denotes the symmetrization with respect to the indices $s_1, \ldots , s_u$.} covariant derivatives}
\begin{eqnarray}
\left({\nabla}^{(u)}R\right)_{i j k l s_1 \ldots s_u} & := & {\nabla}_{(s_1} {\nabla}_{s_2} \ldots {\nabla}_{s_u)} R_{i j k l}\;=\;R_{i j k l\,;\,(s_1 \ldots s_u)}
\end{eqnarray}
are generated by special Young symmetrizers.
\begin{Thm} \label{thm2.3}
Consider the Levi-Civita connection $\nabla$ of a pseudo-Riemannian metric $g$.
For $u \ge 0$ the Riemann tensor and its symmetrized covariant derivatives
${\nabla}^{(u)} R$ fulfil
\begin{eqnarray}
e_t^{\ast} {\nabla}^{(u)} R & = & {\nabla}^{(u)} R
\end{eqnarray}
where $e_t := y_t (u+1)/(2\cdot (u+3)!)$ is an idempotent which is formed from the Young symmetrizer $y_t$ of the standard tableau
\begin{eqnarray} \label{equ2.12}%
t & = &
\begin{array}{|c|c|c|cc|c|}
\hline
1 & 3 & 5 & \ldots & \ldots & (u+4) \\
\hline
2 & 4 & \multicolumn{4}{l}{\hspace{3cm}} \\
\cline{1-2}
\end{array} \,.
\end{eqnarray}
\end{Thm}
A proof of this result of \cite{full4} can be found in \cite[Sec.6]{fie12}, too. The proof needs only the symmetry properties (\ref{equ1.1}) or (\ref{equ1.6}), the identities Bianchi I and Bianchi II and the symmetry with respect to $s_1, \ldots , s_u$. Thus Theorem \ref{thm2.3} is a statement about algebraic curvature tensors and algebraic covariant derivative curvature tensors.
We can specify this in the following way:
\begin{Def}
A tensor $\fR^{(u)} \in {\calT}_{4+u} V$, $u\ge 0$, is called a {\itshape symmetric algebraic covariant derivative curvature tensor of order $u$} iff $\fR^{(u)}(w,y,z,x,a_1,\ldots,a_u)$ is symmetric with respect to $a_1,\ldots,a_u$ and fulfills
\begin{eqnarray}
\fR^{(u)}(w,x,y,z,a_1,\ldots,a_u) & = & - \fR^{(u)}(w,x,z,y,a_1,\ldots,a_u) \\
 & = & \fR^{(u)}(y,z,w,x,a_1,\ldots,a_u) \nonumber
\end{eqnarray}
\begin{eqnarray}
\hspace*{1cm}0 & = & \fR^{(u)}(w,x,y,z,a_1,\ldots,a_u) + \fR^{(u)}(w,y,z,x,a_1,\ldots,a_u) + \\
 & & \fR^{(u)}(w,z,x,y,a_1,\ldots,a_u) \nonumber \\
0 & = & \fR^{(u)}(w,x,y,z,a_1,a_2,\ldots,a_u) + \fR^{(u)}(w,x,z,a_1,y,a_2,\ldots,a_u) + \\
 & & \fR^{(u)}(w,x,a_1,y,z,a_2,\ldots,a_u) \nonumber
\end{eqnarray}
for all $a_1,\ldots, a_u, w, x, y, z\in V$.
\end{Def}
Now symmetric algebraic covariant derivative curvature tensors can be characterized by means of the Young symmetrizer of the tableau (\ref{equ2.12}).
\begin{Prop} \label{prop2.8a}%
A tensor $T \in {\calT}_{4+u} V$, $u\ge 0$, is a symmetric algebraic covariant derivative curvature tensor of order $u$ iff $T$ satisfies
\begin{eqnarray}
e_t^{\ast} T & = & T
\end{eqnarray}
where $e_t$ is the idempotent from Theorem {\rm \ref{thm2.3}}.
\end{Prop}
A proof of Proposition \ref{prop2.8a} is given in the proof of \cite[Prop.6.1]{fie12}. If we consider now the values $u = 0$ and $u = 1$, we obtain
\begin{Cor} \label{cor2.4}%
A tensor $T \in {\calT}_4 V$ $[\tilde{T} \in {\calT}_5]$ is an algebraic {\rm [}covariant derivative{\rm ]} curvature tensor iff $T$ $[\tilde{T}]$ satisfies
\begin{eqnarray}
y_t^{\ast} T\;=\;12\,T & \;\;\;,\;\;\; &
\left[\;y_{t'}^{\ast} \tilde{T}\;=\;24\,\tilde{T}\;\right]
\end{eqnarray}
where $y_t$ $[ y_{t'} ]$ is the Young symmetrizer of the standard tableau
\begin{eqnarray} \label{equ2.16}%
t\;=\;
\begin{array}{|c|c|}
\hline
1 & 3 \\
\hline
2 & 4 \\
\hline
\end{array}
& \;\;\;,\;\;\; &
\left[\;\;t'\;=\;
\begin{array}{|c|c|c|c}
\cline{1-3}
1 & 3 & 5 & \\
\cline{1-3}
2 & 4 &\multicolumn{2}{c}{\;\;\;} \\
\cline{1-2}
\end{array}
\;\;\right]\,.
\end{eqnarray}
\end{Cor}
\vspace{0.5cm}

In the situation considered in Theorem \ref{thm1.10} the symmetry class of the tensors $U$ is not unique. The $(2\,1)$-equivalence class of minimal right ideals $\frakr\subset\bbK[\calS_3]$ which we use to define symmetry classes for the $U$ is an infinite set. The set of generating idempotents for these right ideals $\frakr$ is infinite, too. In \cite{fie03b} we used {\itshape discrete Fourier transforms} to determine a family of primitive generating idempotents of the above minimal right ideals $\frakr\subset\bbK[\calS_3]$.

We denote by $\bbK^{d\times d}$ the set of all $d\times d$-matrices of elements of $\bbK$.
\begin{Def}
A {\it discrete Fourier transform}\footnote{See M. Clausen and U. Baum \cite{clausbaum1,clausbaum2} for details about fast discrete Fourier transforms.} for $\calS_r$ is an isomorphism
\begin{eqnarray}
D : \; \bbK [\calS_r] & \rightarrow &
\bigotimes_{\lambda \vdash r} {\bbK}^{d_{\lambda} \times d_{\lambda}} \label{equ3.12}%
\end{eqnarray}
according to Wedderburn's theorem which maps the group ring $\bbK [\calS_r]$ onto an outer direct product $\bigotimes_{\lambda \vdash r} {\bbK}^{d_{\lambda} \times d_{\lambda}}$ of full matrix rings ${\bbK}^{d_{\lambda} \times d_{\lambda}}$. We denote by $D_{\lambda}$ the {\it natural projections}
$D_{\lambda} : \bbK [\calS_r] \rightarrow
{\bbK}^{d_{\lambda} \times d_{\lambda}}$.
\end{Def}
A discrete Fourier transform maps every group ring element $a\in\bbK[\calS_r]$
to a block diagonal matrix
\begin{eqnarray} \label{equ3.13}%
D :\;\;a\;=\;\sum_{p\in\calS_r}\,a(p)\,p & \mapsto &
\left(
\begin{array}{cccc}
A_{{\lambda}_1} & & & 0 \\
 & A_{{\lambda}_2} & & \\
 & & \ddots & \\
0 & & & A_{{\lambda}_k} \\
\end{array}
\right)\,.
\end{eqnarray}
The matrices $A_{\lambda}$ are numbered by the partitions $\lambda\vdash r$. The dimension $d_{\lambda}$ of every matrix $A_{\lambda}\in{\bbK}^{d_{\lambda} \times d_{\lambda}}$ can be calculated from the Young frame belonging to $\lambda\vdash r$ by means of the {\itshape hook length formula}. For $r = 3$ we have
\begin{center}
\begin{tabular}{c|ccc}
$\lambda$ & $(3)$ & $(2\,1)$ & $(1^3)$ \\
\hline
$d_{\lambda}$ & 1 & 2 & 1 \\
\end{tabular}
\end{center}
The inverse discrete Fourier transform is given by
\begin{Prop}\footnote{See M. Clausen and U. Baum \cite[p.81]{clausbaum1}.}
If $D : \bbK [\calS_r]\rightarrow
\bigotimes_{\lambda \vdash r} {\bbK}^{d_{\lambda} \times d_{\lambda}}$
is a discrete Fourier transform for $\bbK[\calS_r]$, then we have for every
$a\in\bbK[\calS_r]$
\begin{eqnarray}
\forall\,p\in\calS_r:\;\;\;a(p) & = & \frac{1}{r!}\,\sum_{\lambda\vdash r}\,
d_{\lambda}\,\mathrm{trace}\left\{D_{\lambda}(p^{-1})\cdot D_{\lambda}(a)\right\} \label{equ3.14}\\
 & = & \frac{1}{r!}\,\sum_{\lambda\vdash r}\,
d_{\lambda}\,\mathrm{trace}\left\{D_{\lambda}(p^{-1})\cdot A_{\lambda}\right\}\,. \nonumber
\end{eqnarray}
\end{Prop}
In our considerations we are interested in the matrix ring $\bbK^{2\times 2}$ which corresponds to the $(2\,1)$-equivalence class of minimal right ideals $\frakr\subset\bbK[\calS_3]$. In \cite{fie03b} we proved
\begin{Prop}
Every minimal right ideal $\frakr\subset\bbK^{2\times 2}$ is generated by exactly one of the following (primitive) idempotents
\begin{eqnarray} \label{equ3.15}%
Y\;:=\;\left(
\begin{array}{cc}
0 & 0 \\
0 & 1 \\
\end{array}
\right) 
 & \;\;\;or\;\;\; &
X_{\nu}\;:=\;\left(
\begin{array}{cc}
1 & 0 \\
\nu & 0 \\
\end{array}
\right)\;\;\;,\;\;\;\nu\in\bbK\,.
\end{eqnarray}
\end{Prop}
Using an inverse discrete Fourier transform we can
determine the primitive idempotents $\eta, {\xi}_{\nu}\in\bbK[\calS_3]$ which correspond to $Y$, $X_{\nu}$ in (\ref{equ3.15}). We calculated these idempotents by means of the tool
\verb|InvFourierTransform| of the Mathematica package {\sf PERMS} \cite{fie10} (see also \cite[Appendix B]{fie16}.) This calculation can be carried out also by the program package {\sf SYMMETRICA} \cite{kerbkohn2,kerbkohnlas}.
\begin{Prop}
Let us use Young's natural representation\footnote{Three discrete Fourier transforms (\ref{equ3.12}) are known for symmetric groups $\calS_r$: (1) {\itshape Young's natural representation}, (2) {\itshape Young's seminormal representation} and (3) {\itshape Young's orthogonal representation}. See \cite{boerner,boerner2,kerber,clausbaum1}. A short description of (1) and (2) can be found in \cite[Sec.I.5.2]{fie16}. All three discrete Fourier transforms are implemented in the program package {\sf SYMMETRICA} \cite{kerbkohn2,kerbkohnlas}. {\sf PERMS} \cite{fie10}  uses the natural representation. The fast {DFT}-algorithm of M. Clausen and U. Baum \cite{clausbaum1,clausbaum2} is based on the seminormal representation.} of $\calS_3$ as discrete Fourier transform. If we apply the Fourier inversion formula {\rm (\ref{equ3.14})} to a $4\times 4$-block matrix
\begin{eqnarray} \label{equ3.19}%
\left(
\begin{array}{ccc}
A_{(3)} & & 0 \\
 & A_{(2\,1)} & \\
0 & & A_{(1^3)} \\
\end{array}
\right)
 & = &
\left(
\begin{array}{ccc}
0 & 0 & 0 \\
0 & A & 0 \\
0 & 0 & 0 \\
\end{array}
\right)
\end{eqnarray}
where $A$ is equal to $X_{\nu}$ or $Y$ in {\rm (\ref{equ3.15})}, then we obtain the following idempotents of $\bbK[\calS_3]$:
\begin{eqnarray}
X_{\nu}\;\;\;\Rightarrow\;\;\;{\xi}_{\nu} & := &
{\textstyle\frac{1}{3}}\,\{[1,2,3] + \nu [1,3,2] + (1-\nu)[2,1,3] \label{equ3.20}\\
 & & - \nu [2,3,1] + (-1+\nu)[3,1,2] - [3,2,1]\} \nonumber \\
Y\;\;\;\Rightarrow\;\;\;\eta & := &
{\textstyle\frac{1}{3}}\,\{[1,2,3] - [2,1,3] - [2,3,1] + [3,2,1]\}\,. \label{equ3.21}
\end{eqnarray}
\end{Prop}
\begin{Rem}
It is interesting to clear up the connection of the idempotents ${\xi}_{\nu}$ and $\eta$ with Young symmetrizers.
A simple calculation shows that
\begin{eqnarray}
{\xi}_0\;=\;{\textstyle\frac{1}{3}}\,y_{t_1}
 & , &
\eta\;=\;{\textstyle\frac{1}{3}}\,y_{t_2}
\end{eqnarray}
where $y_{t_1}$ and $y_{t_2}$ are the Young symmetrizers of the tableaux
\begin{eqnarray*}
t_1\;=\;
\begin{array}{|c|c|c}
\cline{1-2}
1 & 2 & \\
\cline{1-2}
3 & \multicolumn{2}{c}{\;\;\;} \\
\cline{1-1}
\end{array}
 & , &
t_2\;=\,
\begin{array}{|c|c|c}
\cline{1-2}
1 & 3 & \\
\cline{1-2}
2 & \multicolumn{2}{c}{\;\;\;} \\
\cline{1-1}
\end{array}
\,.
\end{eqnarray*}
\end{Rem}
\begin{Rem}
The proof ot Theorem \ref{thm1.10} is based on the following idea. To treat expressions $y_{t'}^{\ast}(S\otimes U)$ and $y_{t'}^{\ast}(A\otimes U)$ we form the following elements of $\bbK[\calS_5]$:
\begin{eqnarray}
{\sigma}_{\nu , \epsilon} & := & y_{t'}^{\ast}\cdot {\zeta}_{\epsilon}'\cdot
{\xi}_{\nu}''
\;\;\;,\;\;\;
{\rho}_{\epsilon}\;:=\;y_{t'}^{\ast}\cdot {\zeta}_{\epsilon}'\cdot
{\eta}'' \label{equ3.24}\\
{\zeta}_{\epsilon}' & := & \id + \epsilon\,(1\,2)\;\;\;,\;\;\;
\epsilon\in\{1,-1\}\\
{\xi}_{\nu} & \mapsto & {\xi}_{\nu}''\in\bbK[\calS_5]
\;\;\;,\;\;\;
{\eta}\;\mapsto\;{\eta}''\in\bbK[\calS_5]\,. \label{equ3.26}
\end{eqnarray}
Formula (\ref{equ3.26}) denotes the embedding of ${\xi}_{\nu}\,,\,\eta\in\bbK[\calS_3]$ into $\bbK[\calS_5]$ which is induced by the mapping
$\calS_3\rightarrow\calS_5\;,\;
[i_1,i_2,i_3]\mapsto [1,2,i_1 + 2,i_2 + 2,i_3 + 2]$.

For expressions $y_{t'}^{\ast}(U\otimes S)$ and $y_{t'}^{\ast}(U\otimes A)$ we consider the $\bbK[\calS_5]$-elements
\begin{eqnarray}
{\sigma}_{\nu , \epsilon} & := & y_{t'}^{\ast}\cdot {\zeta}_{\epsilon}''\cdot
{\xi}_{\nu}'
\;\;\;,\;\;\;
{\rho}_{\epsilon}\;:=\;y_{t'}^{\ast}\cdot {\zeta}_{\epsilon}''\cdot
{\eta}' \label{equ3.27}\\
{\zeta}_{\epsilon}'' & := & \id + \epsilon\,(4\,5)\;\;\;,\;\;\;
\epsilon\in\{1,-1\}\\
{\xi}_{\nu} & \mapsto & {\xi}_{\nu}'\in\bbK[\calS_5]
\;\;\;,\;\;\;
{\eta}\;\mapsto\;{\eta}'\in\bbK[\calS_5]\,. \label{equ3.29}
\end{eqnarray}
and the embedding
$\calS_3\rightarrow\calS_5\;,\;
[i_1,i_2,i_3]\mapsto [i_1,i_2,i_3,4,5]$
in (\ref{equ3.29}).

Using the Mathematica package {\sf PERMS} \cite{fie10} we verified in \cite{fie03b} that
\begin{center}
\fbox{
${\rho}_{\epsilon}\not= 0 \;\;\;\;\;\;\mathrm{and}\;\;\;\;\;\;
{\sigma}_{\nu , \epsilon}\not= 0\;\Leftrightarrow\;\nu\not=\frac{1}{2}
$}
\end{center}
in both cases.
If ${\rho}_{\epsilon}\not= 0$, ${\sigma}_{\nu , \epsilon}\not= 0$ then the
minimal right ideals $y_{t'}^{\ast}\cdot\bbK[\calS_5]$, ${\rho}_{\epsilon}\cdot\bbK[\calS_5]$ and  ${\sigma}_{\nu , \epsilon}\cdot\bbK[\calS_5]$
coincide, i.e. the symmetry operators ${\rho}_{\epsilon}$, ${\sigma}_{\nu , \epsilon}$ can be used to define the symmetry class of algebraic covariant derivative curvature tensors.
A tensor $T\in\calT_5 V$ is an algebraic covariant derivative curvature tensor iff a tensor $T'\in\calT_5 V$ exists such that
$T = {\rho}_{\epsilon} T'$ or $T = {\sigma}_{\nu , \epsilon} T'$ (if $\nu\not=\frac{1}{2}$).

On the basis of this fact the statement of Theorem \ref{thm1.10} can be proved (see \cite{fie03b}). If $\nu = \frac{1}{2}$, then ${\xi}_{\nu}$ generates the right ideal $\frakr_0 = {\xi}_{\nu}\cdot\bbK[\calS_3]$.  $\frakr_0$ was excluded in Theorem \ref{thm1.10} because ${\sigma}_{1/2 , \epsilon} = 0$.
\end{Rem}\vspace{20pt}

\section{Procedures for the construction of short coordinate representations of the tensors (\ref{equ1.22})} \label{sec4}

In this section we begin to construct short coordinate representations of the tensors (\ref{equ1.22}). We use a common symbol $W_{i j}$ for the coordinates $A_{i j}$ and $S_{i j}$ of the tensors $A, S\in\calT_2 V$. Then the relations
\begin{eqnarray} \label{equ4.1}
S_{i j}\;=\;S_{j i} & , & A_{i j}\;=\;- A_{j i}\,.
\end{eqnarray}
can be written as
\begin{eqnarray}
W_{i j} & = & \epsilon\,W_{j i}\;\;\;,\;\;\;\epsilon\;=\;\left\{
\begin{array}{rl}
1 & {\rm if}\;\;W\;\;{\rm symmetric}\\
-1 & {\rm if}\;\;W\;\;\mbox{{\rm skew-symmetric}}\,.\\
\end{array}
\right.
\end{eqnarray}

If we apply the symmetry operator $\frac{1}{24}y_{t'}^{\ast}$ to tensors $S\otimes U$, $A\otimes U$, $U\otimes S$ or $U\otimes A$, then we obtain long polynomials
\begin{eqnarray} \label{equ4.3a}
\mathfrak{P}_{i_1\ldots i_5} & = & \sum_{p\in\calS_5}\,c_p\,U_{i_{p(1)} i_{p(2)} i_{p(3)}} W_{i_{p(4)} i_{p(5)}}\;\;\;\;,\;\;\;\;c_p\in\bbK
\end{eqnarray}
in the coordinates of $A$, $S$ and $U$. A {\it first reduction} of the number of summands of (\ref{equ4.3a}) can be determined by means of the relations (\ref{equ4.1}).
The results are the expressions which we present in Appendix A. They are polynomials consisting of 40 summands in general.
In particular they fulfil
\begin{eqnarray} \label{equ3.4}%
\mathfrak{P}_{i_1\ldots i_5} & = &
{\textstyle\frac{1}{24}}(y_{t'}^{\ast}(U\otimes W))_{i j k l r} \;=\;
\epsilon {\textstyle\frac{1}{24}}(y_{t'}^{\ast}(W\otimes U))_{i j k l r}\,,
\end{eqnarray}
i.e. a permutation of $U$ and $W$ causes at most a change of the sign (see the lemma in Appendix A). Equation (\ref{equ3.4}) yields statement (1) of Theorem \ref{thm1.11}.

A {\it second reduction} of the number of summands can be carried out by means of all linear identities which the coordinates of the tensor $U$ satisfy.
The determination of the set of these identities is based on the following method from \cite[Sec. III.4.1, I.1.2]{fie16} (see also \cite{fie18}).

Linear identities for the coordinates $T_{i_1\ldots i_r}$ of a tensor $T\in\calT_r V$ have the form
\begin{eqnarray}
0 & = & \sum_{p\in\calS_r}\,x_p\,T_{i_{p(1)}\ldots i_{p(r)}}\;\;\;,\;\;\;
x_p\in\bbK\,,
\end{eqnarray}
where $x_p$ are given numbers from $\bbK$.

Assume that the symmetry class $\calT_{\frakr}$ of $T\in\calT_r V$ is defined by a right ideal $\frakr\subseteq\bbK[\calS_r]$. Then we know from Proposition \ref{prop2.4} that every group ring element $T_b\in\bbK[\calS_r]$, $b\in V^r$, belongs to the left ideal $\frakl = \frakr^{\ast}$. Moreover, $\frakl$ is the smallest linear subspace of $\bbK[\calS_r]$ which contains all $T_b$ of all $T\in\calT_{\frakr}$ if $\dim V\ge r$ (see Proposition \ref{prop2.5}).

Let us denote by $\frakl^{\bot}$ the orthogonal subspace of $\frakl$, i.e. the space of all linear functionals $\langle x , . \rangle$ on $\bbK[\calS_r]$ that vanish on all elements of $\frakl$. Obviously, every $x\in\frakl^{\bot}$ yields a linear identity for the coordinates of $T$ since we can write
\begin{eqnarray} \label{equ4.3}%
0 & = & \langle x\,,\,T_b\rangle\;=\;\sum_{p\in\calS_r}\,x_p\,T_b(p)\;=\;\sum_{p\in\calS_r}\,x_p\,T_{i_{p(1)}\ldots i_{p(r)}}\,,
\end{eqnarray}
where $x_p := \langle x , p \rangle$, $p\in\calS_r$. (The last step is
correct if the $b$ occurring in (\ref{equ4.3}) is an $r$-tuple of basis 
vectors of $V$.)

Every identity (\ref{equ4.3}) can be used to eliminate 
coordinates of $T$ in a polynomial in tensor coordinates. If $\frakl$ is spanned by all $T_b$ of the tensors 
considered , then ${\frakl}^{\bot}$ contains all linear identities which are possible between coordinates of $T$ (compare Prop. \ref{prop2.5}).

We see immediately
\begin{Prop} \label{prop4.1}%
If a basis
$\{ h_1 , \ldots , h_d \}$ of $\frakl$ is known, then
the coefficients $x_p = \langle x , p \rangle$ of the $x \in {\frakl}^{\bot}$ can be obtained
from the linear equation system
\begin{eqnarray}
\langle x , h_i \rangle \;=\;
\sum_{p \in {\calS}_r} h_i(p)\,x_p & = & 0 \hspace{1cm}(i = 1 , \ldots , d)
\,. \label{equ4.4}%
\end{eqnarray}
\end{Prop}
The determination of a basis $\{ h_1 , \ldots , h_d \}$ of $\frakl$ can be carried out efficiently by means of inverse Fourier transforms since there is a fast construction of bases in $\bigotimes_{\lambda\vdash r}\bbK^{d_{\lambda}\times d_{\lambda}}$. We present here only the construction of bases of minimal left ideals of a matrix ring $\bbK^{d\times d}$. See \cite[Sec.I.1.2]{fie16} for more general cases.
\begin{Def}
Let $a\in\bbK^d$ be a $d$-tuple and $i\in\bbN$ be a natural number with $1\le i\le d$. Then we denote by $C_{i,a}\in\bbK^{d\times d}$ that matrix in which the $i$-th row is equal to $a$ and all other rows are filled with zeros.
\end{Def}
\begin{Prop} \label{prop4.3}%
Let $A\in\bbK^{d\times d}$, $A\not= 0$, be a generating element of a minimal left ideal $\frakl = \bbK^{d\times d}\cdot A$ of $\bbK^{d\times d}$. If $a\not= 0$ is a non-vanishing row of $A$, then the matrix set
\begin{eqnarray}
\calB & := & \left\{\,C_{i,a}\;|\;i = 1,\ldots,d\,\right\}
\end{eqnarray}
is a basis of $\frakl$.
\end{Prop}
Proposition \ref{prop4.3} is a special case of \cite[Prop.I.1.33]{fie16} (compare also \cite[Prop.I.1.35]{fie16} or \cite[Prop.5.1]{fie18}).
Due to Proposition \ref{prop4.3} we can construct bases for Proposition \ref{prop4.1} by the following
\begin{Proc} \label{proc4.4}%
Let $e\in\bbK[\calS_r]$ be a generating idempotent of a $d$-dimensional, minimal right ideal $\frakr\subset\bbK[\calS_r]$ which defines a symmetry class $\calT_{\frakr}$ of tensors from $\calT_r V$. Then we can obtain a basis $\{h_1,\ldots,h_d\}$ for Proposition \ref{prop4.1} by the following steps:
\begin{enumerate}
\item{Calculate the generating idempotent $e^{\ast}$ of the minimal left ideal $\frakl = \frakr^{\ast}$.}
\item{Form the discrete Fourier transform of $e^{\ast}$ which possesses only one non-vanishing matrix block $E$, i.e.
\begin{eqnarray*}
D(e^{\ast}) & = &
\left(
\begin{array}{ccc}
0 & & \\
 & E & \\
 & & 0 \\
\end{array}
\right)\,.
\end{eqnarray*}
}
\item{Search for a row $a\not= 0$ of $E$ and form the basis
$\calB = \{ C_{i , a}\;|\;i = 1,\ldots,d \}$ of the left ideal
$\bbK^{d\times d}\cdot E$.}
\item{Calculate
\begin{eqnarray*}
h_i & = & D^{-1}\,\left(
\begin{array}{ccc}
0 & & \\
 & C_{i , a} & \\
 & & 0 \\
\end{array}\right)\,.
\end{eqnarray*}}
\end{enumerate}
\end{Proc}
\begin{Rem}
In the $\calS_3$ we could also determine a basis of $\frakl = \bbK[\calS_3]\cdot e^{\ast}$ by searching for linearly independent vectors in the set
\begin{eqnarray}
 & & \left\{\,p\cdot e^{\ast}\;|\;p\in\calS_r\,\right\}\,. \label{equ4.6}
\end{eqnarray}
However, this way is not efficient in a large $\calS_r$ since (\ref{equ4.6}) leads to an $(r! \times r!)$-matrix of coefficients. In a large $\calS_r$ Procedure \ref{proc4.4} is better than an investigation of (\ref{equ4.6}).
\end{Rem}
If we know a basis $\calB = \{h_1,\ldots,h_d\}$ of $\frakl =\bbK[\calS_r]\cdot e^{\ast}$, then we can determine solutions of the equation system (\ref{equ4.4}) and identities (\ref{equ4.3}) by the following obvious procedure:
\begin{Proc} \label{proc4.6}%
\begin{enumerate}
\item{Choose a subset $\calP\subseteq{\calS_r}$ of $d = \dim\frakl$ permutations, such that $\{n_p := [h_i(p)]_{i = 1,\ldots,d}\;|\;p\in\calP\}$ is a set of $d$ linearly independent column vectors. Then the equation system (\ref{equ4.4}) can be written in the form
\begin{eqnarray} \label{equ4.7}%
0 & = & \sum_{p\in\calP}\,h_i(p)\,x_p + \sum_{p\in{\calS_r}\setminus\calP}\,h_i(p)\,x_p\hspace{1cm}(i = 1 , \ldots , d)\,.
\end{eqnarray}
}
\item{For every fixed $\bar{p}\in{\calS_r}\setminus\calP$ we can determine a unique solution $x_p^{(\bar{p})}$ of (\ref{equ4.7}) which fulfils
\begin{eqnarray}
x_p^{(\bar{p})} & = & \left\{
\begin{array}{cl}
1 & \mathrm{if}\;\;p = \bar{p}\\
0 & \mathrm{if}\;\;p\in{\calS_r}\setminus (\calP\cup\{\bar{p}\})\\
\end{array}
\right.\,.
\end{eqnarray}
The $x_p^{(\bar{p})}$ with $p\in\calP$ can be calculated from the equation system
\begin{eqnarray}
0 & = & \sum_{p\in\calP}\,h_i(p)\,x_p + h_i(\bar{p})\hspace{1cm}(i = 1 , \ldots , d)\,.
\end{eqnarray}
Obviously, the set of all such solutions $\{x_p^{(\bar{p})}\;|\;p\in\calS_r\}$, $\bar{p}\in{\calS_r}\setminus\calP$ is a basis of the null space of (\ref{equ4.7}).
}
\item{Every above solution $\{x_p^{(\bar{p})}\;|\;p\in\calS_r\}$, $\bar{p}\in{\calS_r}\setminus\calP$ defines a linear identity
\begin{eqnarray} \label{equ4.10a}%
0 & = & \sum_{p\in\calP}\,x_p^{(\bar{p})}\,T_{i_{p(1)}\ldots i_{p(r)}} +
T_{i_{\bar{p}(1)}\ldots i_{\bar{p}(r)}}\hspace{1cm}(\bar{p}\in{\calS_r}\setminus\calP)
\end{eqnarray}
for the coordinates of the tensors $T$ of the symmetry class $\calT_{\frakl^{\ast}}$.
}
\end{enumerate}
\end{Proc}
The identities (\ref{equ4.10a}) have the remarkable property that they depend only on the choice of the set $\calP$.
\begin{Thm}
Let $\calB = \{h_1,\ldots,h_d\}$ and $\tilde{\calB} = \{\tilde{h}_1,\ldots,\tilde{h}_d\}$ be bases of the above $d$-dimensional left ideal $\frakl$, which defines a symmetry class $\calT_{\frakl^{\ast}}$. Furthermore, let $\calP\subseteq\calS_r$ be a subset of $d$ permutations such that
$\{n_p = [h_i(p)]_{i=1,\ldots,d}\;|\;p\in\calP\}$ is a set of $d$ linearly independent column vectors. Then the following holds:
\begin{enumerate}
\item{The set
$\{\tilde{n}_p = [\tilde{h}_i(p)]_{i=1,\ldots,d}\;|\;p\in\calP\}$ belonging to $\tilde{\calB}$ also consists of $d$ linearly independent column vectors.
}
\item{For $\calB$, $\tilde{\calB}$ and $\calP$ the Procedure {\rm \ref{proc4.6}} yields sets of identities {\rm (\ref{equ4.10a})} whose coefficients satisfy
\begin{eqnarray}
\forall\,p\in\calP:\;\;\;x_p^{(\bar{p})} & = & \tilde{x}_p^{(\bar{p})}\,. \label{equ4.11a}
\end{eqnarray}
}
\end{enumerate}
\end{Thm}
\begin{proof}
Ad (1): The vectors $h_i$, $\tilde{h}_i$ of $\calB$, $\tilde{\calB}$ satisfy a relation
\begin{eqnarray}
\tilde{h}_i & = & \sum_{k = 1}^d\,K_{ik}\,h_k \label{equ4.12}
\end{eqnarray}
where $K := [K_{ik}]_{i,k = 1,\ldots,d}$ is a regular $d\times d$-matrix. From (\ref{equ4.12}) we obtain $\tilde{n}_p = K\cdot n_p$ for all $p\in\calP$. Consequently the $\tilde{n}_p$ are linearly independent again.

Ad (2): For all $\bar{p}\in\calS_r\setminus\calP$ the $\tilde{x}_p^{(\bar{p})}$ with $p\in\calP$ can be calculated from the equation system
\begin{eqnarray}
0 & = & \sum_{p\in\calP}\,\tilde{h}_i(p)\,\tilde{x}_p + \tilde{h}_i(\bar{p})\hspace{1cm}(i = 1 , \ldots , d) \label{equ4.13}
\end{eqnarray}
which is equivalent to
\begin{eqnarray}
0 & = & \sum_{k = 1}^d K_{ik}\,\left\{\sum_{p\in\calP}\,h_k(p)\,\tilde{x}_p + h_k(\bar{p})\right\}\hspace{1cm}(i = 1 , \ldots , d)\,. \label{equ4.14}
\end{eqnarray}
Since the $\tilde{n}_p$ are linearly independent column vectors, system (\ref{equ4.13}) has a unique solution. On the other hand, we see from (\ref{equ4.14}) that the $x_p^{(\bar{p})}$ with $p\in\calP$ solve (\ref{equ4.14}) and (\ref{equ4.13}). Thus (\ref{equ4.11a}) is correct.
\end{proof}
Now we consider a finite set of $T^{(1)}$, $T^{(2)}$, \ldots , $T^{(m)}$ of covariant tensors of orders $r_1$, $r_2$, \ldots , $r_m$ and a polynomial
\begin{eqnarray} \label{equ4.17}
\mathfrak{P}(T^{(1)},\ldots , T^{(m)})_{i_1 \ldots i_r} & = &
\sum_{q\in\calS_r}\,c_q\,T^{(1)}_{i_{q(1)}\ldots i_{q(r_1)}}\cdot\ldots\cdot
T^{(m)}_{i_{q(r-r_m+1)}\ldots i_{q(r)}} \\
 & = &
\sum_{q\in\calS_r}\,c_q\,(T^{(1)}\otimes\ldots\otimes
T^{(m)})_{i_{q(1)}\ldots i_{q(r)}} \nonumber
\end{eqnarray}
where $r = r_1 +\ldots + r_m$. We do not require that $c_q\not= 0$ for all $q\in\calS_r$. Thus the summation can run through a subset of $\calS_r$ only.
Furthermore we assume that a {\it lexicographic ordering} is defined for the index names in (\ref{equ4.17}).

If we know identities of type (\ref{equ4.10a}) for a tensor $T^{(k)}$ in (\ref{equ4.17}), then we can transform all coordinates
$\{T^{(k)}_{j_{s(1)}\ldots j_{s(r_k)}}\;|\; s\in\calS_{r_k}\}$ which belong to the same set of index names $\{j_1, \ldots , j_{r_k}\}\subset\{i_1, \ldots  , i_r\}$ into the set $\{ T^{(k)}_{j_{p(1)} \ldots j_{p(r_k)}} \;|\; p\in\calP\}$ of $d = |\calP|$ coordinates. Here we assume that $\{j_1, \ldots , j_{r_k}\}$ is an arrangement of indices whose names form an increasing sequence according to the lexicographic ordering.

The identities (\ref{equ4.10a}) do not depend on the choice of a basis $\calB$, but only on the choice of the set $\calP$. We can carry out the above reduction
\begin{eqnarray} \label{equ4.18}
\{T^{(k)}_{j_{s(1)}\ldots j_{s(r_k)}}\;|\; s\in\calS_{r_k}\}
& \rightarrow &
\{ T^{(k)}_{j_{p(1)} \ldots j_{p(r_k)}} \;|\; p\in\calP\}
\end{eqnarray}
of tensor coordinates for every set $\calP$ which satisfies the assumtions of Procedure \ref{proc4.6}. However, if we apply this reduction process to a polynomial (\ref{equ4.17}) and reduce the set of coordinates of a single tensor $T^{(k)}$ in (\ref{equ4.17}), then it may happen that we obtain different lengths of the transformed $\mathfrak{P}_{i_1\ldots i_r}$ for different sets $\calP$.
Consequently there is the problem to find such a set $\calP$ for which the transformed $\mathfrak{P}_{i_1\ldots i_r}$ has a minimal lenght.

In the treatment of polynomials (\ref{equ4.3a}) we have to carry out transformations (\ref{equ4.10a}) only with respect to the tensor $U$. The above discussion leads to the following
\begin{Proc} \label{proc4.8}
Let $\frakl = \bbK[\calS_3]\cdot e^{\ast}$ be a minimal left ideal describing the symmetry class of the tensor $U$ and $\calB = \{h_1 , h_2\}$ be a basis of $\frakl$ determined by Procedure \ref{proc4.4}. (The hook length formula tells us that
$\dim\frakl = 2$ for every minimal left ideal $\frakl$ belonging to the equivalence class of $(2\,1)\vdash 3$.) Then carry out the following steps for every subset $\calP\subset\calS_3$ with $|\calP| = 2$:
\begin{enumerate}
\item{Check the condition
\begin{eqnarray} \label{equ4.19}
\Delta_{\calP}\;:=\;\det \Big(h_i(p)\Big)_{\twoline{i = 1,2}{p\in\calP}} & \not= & 0\,.
\end{eqnarray}
If (\ref{equ4.19}) is not valid then skip the steps (2) and (3) for the set $\calP$.
}
\item{If $\Delta_{\calP}\not= 0$ then determine identities (\ref{equ4.10a}) for the coordinates of $U$ by means of Procedure \ref{proc4.6} from $\calB$ and $\calP$.}
\item{Carry out a reduction (\ref{equ4.18}) of the coordinates of $U$ in (\ref{equ4.3a}) by means of the identities (\ref{equ4.10a}) from step (2). Determine the number of summands of the reduced polynomial (\ref{equ4.3a}).}
\end{enumerate}
\end{Proc}
\begin{Rem} \label{rem4.9}
It is clear that Procedure \ref{proc4.8} can be generalized.
\begin{enumerate}
\item{If a polynomial (\ref{equ4.17}) with $m$ tensors $T^{(k)}$ is given, then
we could reduce the set of coordinates of every tensor $T^{(k)}$ in (\ref{equ4.17}) by a procedure of the type \ref{proc4.8}.}
\item{An exact transfer of Procedure \ref{proc4.8} to the general case (\ref{equ4.17}) means that we determine identities (\ref{equ4.10a}) for a given set $\calP$ and use these fixed identities for all index sets $\{j_1, \ldots , j_{r_k}\}$ of the considered tensor $T^{(k)}$ in the reduction process (\ref{equ4.18}). We could also go a more general way. For every set of index names $\{j_1, \ldots , j_{r_k}\}$ occuring in coordinates of the considered tensor $T^{(k)}$ we could run through all possible sets $\calP$ and search for such a $\calP$ for which the length of (\ref{equ4.17}) becomes minimal. This modification of Procedure \ref{proc4.8} could yield shorter polynomials (\ref{equ4.17}) than the original Procedure \ref{proc4.8}.}
\end{enumerate}
\end{Rem}
In the present paper we use the Procedure \ref{proc4.8} in the above form since the generalization Remark \ref{rem4.9} (2) leads to a considerable increase of the expenditure of calculation. The results of the above Procedure \ref{proc4.8} are polynomials $\mathfrak{P}(U,W)_{i_1\ldots i_5}$ of type (\ref{equ4.3a}) in which {\it the arrangement of the indices of $U$ is defined by one and the same set $\calP$ for all sets $\{j_1, j_2 , j_3\}$ of index names} belonging to this tensor.\vspace{20pt}

\section{Determination of short coordinate representations of the tensors (\ref{equ1.22})} \label{sec5}

Now we apply the Procedures \ref{proc4.4}, \ref{proc4.8} and the Subprocedure \ref{proc4.6} to the idempotents (\ref{equ3.20}), (\ref{equ3.21}) to obtain linear identities (\ref{equ4.10a}) for the tensors $U$ and to determine the shortest polynomials (\ref{equ4.3a}) which can be constructed by these procedures. We carry out our calculations by means of {\sf PERMS} \cite{fie10}. The {\sf Mathematica} notebooks of the computation can be downloaded from \cite{fie21}.

In all following calculations we make the {\it assumtion} that
\begin{eqnarray} \label{equ4.1a}%
 & & \fbox{$\dim V\;\ge\;3$\,.}
\end{eqnarray}
Then the left ideals $\frakl = {\frakr}^{\ast}$ which belong to the symmetry classes of the tensors $U$ are generated by all $U_b$ of the $U$ and the manifold of solutions of the corresponding equation system (\ref{equ4.4}) defines a complete set of linear identities for the tensors $U$ (see Proposition \ref{prop2.5}).

In our calculations we use the following numbering of the permutations of $\calS_3$.
\begin{eqnarray} \label{equ5.2}
 & &
\begin{array}{|c|c|c|c|c|c|c|}
\hline
i & 1 & 2 & 3 & 4 & 5 & 6 \\
\hline
p_i & [1,2,3] & [1,3,2] & [2,1,3] & [2,3,1] & [3,1,2] & [3,2,1] \\
\hline
\end{array}\,.
\end{eqnarray}
\vspace{10pt}

\subsection{Short formulas with a tensor $U$ defined by $\eta$} \label{sec5.1}
The hook lenght formula tells us that the idempotent $\eta$ given in (\ref{equ3.21}) generates a 2-dimensional right ideal.
If we apply Procedure \ref{proc4.4} to $\eta$, we obtain
\begin{eqnarray*}
{\eta}^{\ast} & = &
{\textstyle\frac{1}{3}}\,\{[1,2,3] - [2,1,3] - [3,1,2] + [3,2,1]\}\\
D({\eta}^{\ast}) & = &
\left(
\begin{array}{ccc}
0 & & \\
 & E & \\
 & & 0 \\
\end{array}
\right)
\;=\;
\frac{1}{3}\,\left(
\begin{array}{cccc}
0 & & & \\
 & - 1 & 2 & \\
 & - 2 & 4 & \\
 & & & 0 \\
\end{array}
\right)\,.
\end{eqnarray*}
The first row $a = \frac{1}{3} (-1\,,\,2)$ of $E$ leads to the basis
\begin{eqnarray*}
C_{1 , a}\;=\;\frac{1}{3}\,
\left(
\begin{array}{rr}
-1 & 2 \\
0 & 0 \\
\end{array}
\right)
& , &
C_{2 , a}\;=\;\frac{1}{3}\,
\left(
\begin{array}{rr}
0 & 0 \\
-1 & 2 \\
\end{array}
\right)
\end{eqnarray*}
of $\bbK^{2\times 2}\cdot E$ and an inverse Fourier transform yields the basis
\begin{eqnarray} \label{equ4.8}%
\hspace*{0.5cm}h_1 & = &
{\textstyle\frac{1}{9}}\,\{
- [1,2,3] + 2 [1,3,2] - [2,1,3] + 2 [2,3,1] - [3,1,2] - [3,2,1]\} \\
\hspace*{0.5cm}h_2 & = &
{\textstyle\frac{1}{9}}\,\{
2 [1,2,3] - [1,3,2] - [2,1,3] - [2,3,1] - [3,1,2] + 2 [3,2,1]\} \nonumber
\end{eqnarray}
of $\frakl = \bbK[\calS_3]\cdot{\eta}^{\ast}$.

The coefficients of the group ring elements (\ref{equ4.8}) form the rows of the coefficient matrix of the system (\ref{equ4.4}). From (\ref{equ4.8}) we obtain the matrix
\begin{eqnarray} \label{equ5.4}
 & &
\frac{1}{9}\,\left(
\begin{array}{rrrrrr}
-1 & 2 & -1 & 2 & -1 & -1 \\
2 & -1 & -1 & -1 & -1 & 2 \\
\end{array}
\right)\,.
\end{eqnarray}
Let us carry out a step of Procedure \ref{proc4.8} for the matrix (\ref{equ5.4}). We start with the set $\calP = \{[1,2,3], [1,3,2]\}$ which is connected with the first two columns of (\ref{equ5.4}). Obviously we have ${\Delta}_{\calP}\not= 0$.
The Gaussian algorithm transforms (\ref{equ5.4}) into 
\begin{eqnarray} \label{equ4.9}%
 & &
\left(
\begin{array}{rrrrrr}
1 & 0 & -1 & 0 & -1 & 1 \\
0 & 1 & -1 & 1 & -1 & 0 \\
\end{array}
\right)\,.
\end{eqnarray}
The null space of (\ref{equ4.9}) has the basis
\begin{eqnarray} \label{equ4.10}%
 & &
\begin{array}{cc}
(-1,0,0,0,0,1) & (1,1,0,0,1,0) \\
(0,-1,0,1,0,0) & (1,1,1,0,0,0) \\
\end{array}\;.
\end{eqnarray}
If we use the vectors (\ref{equ4.10}) to define linear identities (\ref{equ4.10a}) for the coordinates of $U\in\calT_{\frakl^{\ast}}$ we obtain
\begin{eqnarray} \label{equ4.11}%
 & &
\begin{array}{ccccccccccc}
- & U_{ijk} & + &         &   &         &         &         & U_{kji} & = & 0 \\
  & U_{ijk} & + & U_{ikj} & + &         &         & U_{kij} &         & = & 0 \\
  &         & - & U_{ikj} & + &         & U_{jki} &         &         & = & 0 \\
  & U_{ijk} & + & U_{ikj} & + & U_{jik} &         &         &         & = & 0 \\
\end{array}\;.
\end{eqnarray}
Using the identities (\ref{equ4.11}) we can express all coordinates of the tensor $U$ by the coordinates $U_{ijk}$ and $U_{ikj}$ which possess the two lexicographically smallest index arrangements (defined by $\calP$). If we carry out the reduction process (\ref{equ4.18}) for
\begin{eqnarray}
\mathfrak{P}_{i_1\ldots i_5} & := &
{\textstyle\frac{1}{24}}(y_{t'}^{\ast}(U\otimes W))_{i j k l r} \;=\;
\epsilon {\textstyle\frac{1}{24}}(y_{t'}^{\ast}(W\otimes U))_{i j k l r}
\end{eqnarray}
by means of (\ref{equ4.11}), then we obtain reduced polynomials $\mathfrak{P}_{i_1\ldots i_5}^{\rm red}$ of length 12 if $\epsilon = 1$ and of length 20 if $\epsilon = -1$.
\begin{table}[t]
\begin{tabular}{c|c|c|c}
\hline
 & & \multicolumn{2}{c}{length of $\mathfrak{P}_{i_1\ldots i_5}^{\rm red}$}\\
$\calP$ & ${\Delta}_{\calP}$ & $\epsilon = 1$ & $\epsilon = -1$ \\
\hline
 12  &    $\not= 0$ &      12  &     20 \\
 13  &    $\not= 0$ &      16  &     20 \\
 14  &    $\not= 0$ &      12  &     20 \\
 15  &    $\not= 0$ &      16  &     20 \\
 16  &     0        &          & \\
 23  &    $\not= 0$ &      12  &     20 \\
 24  &     0        &          & \\
 25  &    $\not= 0$ &      12  &     20 \\
 26  &    $\not= 0$ &      12  &     20 \\
 34  &    $\not= 0$ &      12  &     20 \\
 35  &     0        &          & \\
 36  &    $\not= 0$ &      16  &     20 \\
 45  &    $\not= 0$ &      12  &     20 \\
 46  &    $\not= 0$ &      12  &     20 \\
 56  &    $\not= 0$ &      16  &     20 \\
\hline
\end{tabular}
\vspace{3mm}
\caption{The lengths of $\mathfrak{P}_{i_1\ldots i_5}^{\rm red}$ for an $U$ from the symmetry class given by $\eta$.}
\end{table}

Table 1 shows the results of the calculations for the remaining 14 sets $\calP$. In the first column a set $\calP$ is denoted by two numbers $mn$ if $\calP$ contains the $m$-th and the $n$-th permutation from (\ref{equ5.2}). We see that different lengths of $\mathfrak{P}_{i_1\ldots i_5}^{\rm red}$ occur for symmetric tensors $W$ whereas an alternating tensor $W$ leads only to a length 20. In Appendix B we present the formulas for $\mathfrak{P}_{i_1\ldots i_5}^{\rm red}$ under the conditions of the first row of Table 1 (i.e. $\calP = \{[1,2,3]\,,\,[1,3,2]\}$ and $\epsilon = \pm 1$). The computer calculations for Table 1 and Appendix B can be found in the {\sf Mathematica} notebooks
\verb|part12b.nb|, \ldots , \verb|part56b.nb| in \cite{fie21}.\vspace{10pt}

\subsection{Short formulas with a tensor $U$ defined by ${\xi}_{\nu}$} \label{sec5.2}
Now we search for short formulas $\mathfrak{P}_{i_1\ldots i_5}^{\rm red}$ which can be formed by means of a tensor $U$ from the symmetry class defined by the idempotent ${\xi}_{\nu}$ (see (\ref{equ3.20})).

Again, the hook length formula tells us that the idempotent ${\xi}_{\nu}$ generates a 2-dimensional right ideal of $\bbK[\calS_3]$ for all $\nu\in\bbK$.
If we start Procedure \ref{proc4.4} for ${\xi}_{\nu}$, we obtain
\begin{eqnarray*}
{\xi}_{\nu}^{\ast} & = &
{\textstyle\frac{1}{3}}\,\{[1,2,3] + \nu\,[1,3,2] + (1-\nu)[2,1,3]\\
 & & + (-1+\nu)[2,3,1] - \nu\,[3,1,2] - [3,2,1]\}\\
D({\xi}_{\nu}^{\ast}) & = &
\left(
\begin{array}{ccc}
0 & & \\
 & E & \\
 & & 0 \\
\end{array}
\right)
\;=\;
\frac{1}{3}\,\left(
\begin{array}{cccc}
0 & & & \\
 & 4 - 2\nu & -2 + 4\nu & \\
 & 2 - \nu & -1 + 2\nu & \\
 & & & 0 \\
\end{array}
\right)\,.
\end{eqnarray*}
It is easy to see that the first row of $E$ is unequal to zero for all $\nu\in\bbK$. The condition $4 - 2\nu = 0$ leads to $\nu = 2$. However, we obtain $-2 + 4\nu = 6$ for $\nu = 2$. Thus the first row $a$ of $E$ can be used for the construction of a basis $\calB$ of the left ideal $\bbK^{2\times 2}\cdot E$ according step (3) of Procedure \ref{proc4.4}. We obtain
\begin{eqnarray*}
C_{1 , a}\;=\;\frac{1}{3}\,
\left(
\begin{array}{rr}
4 - 2\nu & -2 + 4\nu \\
0 & 0 \\
\end{array}
\right)
& , &
C_{2 , a}\;=\;\frac{1}{3}\,
\left(
\begin{array}{rr}
0 & 0 \\
4 - 2\nu & -2 + 4\nu \\
\end{array}
\right)\,.
\end{eqnarray*}
Now an inverse Fourier transform yields the basis
\begin{eqnarray} \label{equ5.9}%
\hspace*{0.5cm}h_1 & = &
{\textstyle\frac{1}{9}}\,\{
(4 - 2\nu)[1,2,3] + (-2 + 4\nu)[1,3,2] + (4 - 2\nu)[2,1,3] \nonumber\\
 & & + (-2 + 4\nu)[2,3,1] - (2 + 2\nu)[3,1,2] - (2 + 2\nu)[3,2,1]\} \\
\hspace*{0.5cm}h_2 & = &
{\textstyle\frac{1}{9}}\,\{
(-2 + 4\nu)[1,2,3] + (4 - 2\nu)[1,3,2] - (2 + 2\nu)[2,1,3] \nonumber \\
 & &
 - (2 + 2\nu)[2,3,1] + (4 - 2\nu)[3,1,2] + (-2 + 4\nu)[3,2,1]\} \nonumber
\end{eqnarray}
of $\frakl = \bbK[\calS_3]\cdot{\eta}^{\ast}$.

Again, the coefficients of the group ring elements (\ref{equ5.9}) form the rows of the coefficient matrix of the system (\ref{equ4.4}). From (\ref{equ5.9}) we obtain the matrix
\begin{eqnarray} \label{equ5.10}
 & &
\frac{1}{9}\,\left(
\begin{array}{rrrrrr}
4 - 2\nu & -2 + 4\nu & 4 - 2\nu & -2 + 4\nu & -2 - 2\nu & -2 - 2\nu \\
-2 + 4\nu & 4 - 2\nu & -2 - 2\nu & -2 - 2\nu & 4 - 2\nu & -2 + 4\nu \\
\end{array}
\right)\,.
\end{eqnarray}
Since all elements of the matrix (\ref{equ5.10}) depend linearly on $\nu$, every determinant ${\Delta}_{\calP}$ is a quadratic polynomial in $\nu$. Our computer calculations show that ${\Delta}_{\calP}\not\equiv 0$ for all 15 sets $\calP\subset\calS_3$ with $|\calP| = 2$.  Consequently, there are at most two roots ${\nu}_1 , {\nu}_2\in\bbK$ for a fixed $\calP$ such that ${\Delta}_{\calP}({\nu}_i) = 0$.\vspace{5pt}

\subsubsection{The generic case} \label{sec5.2.1}
Let us first investigate the {\it generic case}. For every $\calP$ we restrict us to such values $\nu\in\bbK$ for which ${\Delta}_{\calP}(\nu)\not= 0$. Then we can carry out the Procedures \ref{proc4.6} and \ref{proc4.8}. Obviously the coefficients of the identities (\ref{equ4.10a}) and of the reduced polynomials ${\mathfrak{P}}_{i_1\ldots i_5}^{\rm red}$ will then be rational functions in $\nu$. More precisely, for every set $\calP\subset\calS_3$ with $|\calP| = 2$ there are polynomials
\[
P_q^{\calP}(\nu)\;,\;Q_q^{\calP}(\nu)\in\bbK[\nu]\;\;\;,\;\;q\in\calS_5\,,
\]
such that the result of Procedure \ref{proc4.8} belonging to $\calP$ can be written as
\begin{eqnarray} \label{equ5.11}
\mathfrak{P}_{i_1\ldots i_5}^{\rm red} & = & \sum_{q\in\calS_5}\,\frac{P_q^{\calP}(\nu)}{Q_q^{\calP}(\nu)}\,U_{i_{q(1)} i_{q(2)} i_{q(3)}} W_{i_{q(4)} i_{q(5)}}\,.
\end{eqnarray}
\begin{Rem} \label{rem5.1a}%
Note that the calculation of (\ref{equ5.11}) is possible for all values $\nu\in\bbK$ which fulfill ${\Delta}_{\calP}(\nu)\not= 0$. In particular the polynomials $Q_q^{\calP}(\nu)$ in (\ref{equ5.11}) satisfy $Q_q^{\calP}(\nu)\not= 0$ for all $\nu$ with ${\Delta}_{\calP}(\nu)\not= 0$. If we would use {\it Cramer's rule} to solve the equation system (\ref{equ4.7}) and to determine the linear identities (\ref{equ4.10a}), we could write $Q_q^{\calP}(\nu) = {\Delta}_{\calP}(\nu)$ for all $Q_q^{\calP}(\nu)$. However, we give up a condition such as $Q_q^{\calP}(\nu) = {\Delta}_{\calP}(\nu)$ to allow a reduction of the fractions $P_q^{\calP}(\nu)/Q_q^{\calP}(\nu)$.
\end{Rem}

Let us carry out Procedure \ref{proc4.8} for the set $\calP = \{[1,2,3]\,,\,[1,3,2]\}$ which characterizes the first two colums of (\ref{equ5.10}). In this case we get the determinant
\begin{eqnarray} \label{equ5.12}
{\Delta}_{\calP}(\nu) & = & {\textstyle\frac{4}{27}}\,(1 - \nu)(1 + \nu)\,.
\end{eqnarray}
(\ref{equ5.12}) has the roots ${\nu}_1 = 1$ and ${\nu}_2 = -1$. For $\nu\not\in\{1\,,\,-1\}$ Procedure \ref{proc4.6} yields now the identities
\begin{eqnarray} \label{equ5.13}%
 & &
\begin{array}{ccccccccccc}
- & \frac{{\nu}^2 - \nu + 1}{{\nu}^2 - 1}\,U_{ijk} & + & \frac{2\nu - 1}{{\nu}^2 - 1}\,U_{ikj}        & +  &         &         &         & U_{kji} & = & 0 \\
  & \frac{{\nu}^2 - 2\nu}{{\nu}^2 - 1}\,U_{ijk} & + & \frac{{\nu}^2 - \nu + 1}{{\nu}^2 - 1}\,U_{ikj} & + &         &         & U_{kij} &         & = & 0 \\
  & \frac{2\nu - 1}{{\nu}^2 - 1}\,U_{ijk}        & - & \frac{{\nu}^2 - \nu + 1}{{\nu}^2 - 1}\,U_{ikj} & + &         & U_{jki} &         &         & = & 0 \\
  & \frac{{\nu}^2 - \nu + 1}{{\nu}^2 - 1}\,U_{ijk} & + & \frac{{\nu}^2 - 2\nu}{{\nu}^2 - 1}\,U_{ikj} & + & U_{jik} &         &         &         & = & 0 \\
\end{array}\;.
\end{eqnarray}
If we apply the identities (\ref{equ5.13}) to reduce the polynomials $\mathfrak{P}_{i_1\ldots i_5}$ in Appendix A, we obtain reduced polynomials $\mathfrak{P}_{i_1\ldots i_5}^{\rm red}$ of length 16 if $W$ is symmetric and of length 20 if $W$ is skew-symmetric. (These numbers are only correct if we consider the generic case in which all polynomials $P_q^{\calP}(\nu)$ in (\ref{equ5.11}) do not vanish.) Explicite formulas for the reduced polynomials $\mathfrak{P}_{i_1\ldots i_5}^{\rm red}$ are presented in Appendix C. The computer calculations can be found in the {\sf Mathematica} notebooks \cite[proc44.nb,part12a.nb]{fie21}.
\begin{table}[t]
\begin{tabular}{c|c|c|c}
\hline
 & & \multicolumn{2}{c}{length of $\mathfrak{P}_{i_1\ldots i_5}^{\rm red}$}\\
$\calP$ & roots of ${\Delta}_{\calP}(\nu)$ & $\epsilon = 1$ & $\epsilon = -1$ \\
\hline
 12  &    $1\,,\,-1$                             &      16  &     20 \\
 13  &    $0\,,\,2$                              &      16  &     20 \\
 14  &    $e^{\imath\pi/3}\,,\,e^{-\imath\pi/3}$ &      16  &     20 \\
 15  &    $e^{\imath\pi/3}\,,\,e^{-\imath\pi/3}$ &      16  &     20 \\
 16  &    $1/2$                                  &      16  &     20 \\
 23  &    $e^{\imath\pi/3}\,,\,e^{-\imath\pi/3}$ &      16  &     20 \\
 24  &    $1/2$                                  &      16  &     20 \\
 25  &    $0\,,\,2$                              &      16  &     20 \\
 26  &    $e^{\imath\pi/3}\,,\,e^{-\imath\pi/3}$ &      16  &     20 \\
 34  &    $1\,,\,-1$                             &      16  &     20 \\
 35  &    $1/2$                                  &      16  &     20 \\
 36  &    $e^{\imath\pi/3}\,,\,e^{-\imath\pi/3}$ &      16  &     20 \\
 45  &    $e^{\imath\pi/3}\,,\,e^{-\imath\pi/3}$ &      16  &     20 \\
 46  &    $0\,,\,2$                              &      16  &     20 \\
 56  &    $1\,,\,-1$                             &      16  &     20 \\
\hline
\end{tabular}
\vspace{3mm}
\caption{The lengths of $\mathfrak{P}_{i_1\ldots i_5}^{\rm red}$ for an $U$ from the symmetry class given by ${\xi}_{\nu}$ in the generic case ${\Delta}_{\calP}(\nu)\not= 0$ and $P_q^{\calP}(\nu)\not= 0$ for all $q\in\calS_5$.}
\end{table}

When we repeat this computation for the remaining 14 sets $\calP$, we obtain the results given in Table 2. The second column of Table 2 contains the roots of ${\Delta}_{\calP}(\nu)$ which we have to exclude. The computer calculations for Table 2 can be found in the {\sf Mathematica} notebooks \verb|part12a.nb|, \ldots , \verb|part56a.nb| in \cite{fie21}.
\begin{Rem} \label{rem5.1}%
Note that some ${\Delta}_{\calP}(\nu)$ possess only the root $\nu = \frac{1}{2}$ which is the critical $\nu$-value for the construction of algebraic covariant derivative curvature tensors. Since all $\nu\not=\frac{1}{2}$ are allowed values in (\ref{equ5.11}) for such sets $\calP$, the $\mathfrak{P}_{i_1\ldots i_5}^{\rm red}$ of such $\calP$ are formulas which represent algebraic covariant derivative curvature tensors for every $\nu\not=\frac{1}{2}$. As an example, we present in Appendix C formulas of $\mathfrak{P}_{i_1\ldots i_5}^{\rm red}$ for the set $\calP$ denoted by ''16'' in Table 2.
\end{Rem}
\begin{Rem}
We see that the lengths of $\mathfrak{P}_{i_1\ldots i_5}^{\rm red}$ are independent on $\calP$ in the generic case. This property is not self-evident. The following considerations show that our calculation could just as well result in contrary findings.

We denote by $I := \{i,j,k,l,r\}$ the set of index names used for the polynomials of tensor coordinates in the Appendices and by '$\prec$' the lexicographic order for these index names. Let ${\mathfrak I}$ be the set of all pairs $(j_1,j_2)\in I\times I$ with $j_1\prec j_2$. For every $(j_1,j_2)\in{\mathfrak I}$ we denote by $i_1,i_2,i_3$ the elements of $I\setminus\{j_1,j_2\}$ where we assume $i_1\prec i_2\prec i_3$.

Using these notations we can write the expression for ${\textstyle\frac{1}{24}}(y_{t'}^{\ast}(U\otimes W))_{i j k l r}$ given in Appendix A in the form
\begin{eqnarray} \label{equ5.14}%
\hspace*{1cm}{\mathfrak P}_{ijklr}\;:=\;{\textstyle\frac{1}{24}}(y_{t'}^{\ast}(U\otimes W))_{i j k l r}
& = &
\sum_{(j_1,j_2)\in{\mathfrak I}}\,\sum_{s\in\calS_3}\,A_{\epsilon , s}^{(j_1,j_2)}\,U_{i_{s(1)} i_{s(2)} i_{s(3)}} W_{j_1 j_2}\,,
\end{eqnarray}
where $A_{\epsilon , s}^{(j_1,j_2)}\in\bbK$ are constant coefficients. (If the summation with respect to $s$ runs only through a subset of $\calS_3$, then some of these coefficients vanish).

If a set $\calP\subset\calS_3$ of two permutations is given, then we can express the linear identities (\ref{equ4.10a}) for $U$ by
\begin{eqnarray} \label{equ5.15}%
U_{i_{s(1)} i_{s(2)} i_{s(3)}} & = &
\sum_{p\in\calP}\,B_{s p}^{\calP}(\nu)\,U_{i_{p(1)} i_{p(2)} i_{p(3)}}
\;\;\;,\;\;\;s\in\calS_3\setminus\calP\,.
\end{eqnarray}
Here $B_{s p}^{\calP}(\nu)$ are rational functions of $\nu$.

When we set (\ref{equ5.15}) into (\ref{equ5.14}), we obtain
\begin{eqnarray*}
{\mathfrak P}_{ijklr}^{\rm red}
& = &
\sum_{(j_1,j_2)\in{\mathfrak I}}\,\sum_{p\in\calP}\,A_{\epsilon , p}^{(j_1,j_2)}\,U_{i_{p(1)} i_{p(2)} i_{p(3)}} W_{j_1 j_2} + \\
 & &
\sum_{(j_1,j_2)\in{\mathfrak I}}\,\sum_{s\in\calS_3\setminus\calP}\,\sum_{p'\in\calP}\,A_{\epsilon , s}^{(j_1,j_2)}\,B_{s p'}^{\calP}(\nu)\,U_{i_{p'(1)} i_{p'(2)} i_{p'(3)}} W_{j_1 j_2} \\
 & = &
\sum_{(j_1,j_2)\in{\mathfrak I}}\,\sum_{p\in\calP}\,\left\{A_{\epsilon , p}^{(j_1,j_2)} + \sum_{s\in\calS_3\setminus\calP}\,A_{\epsilon , s}^{(j_1,j_2)}\,B_{s p}^{\calP}(\nu)\right\}\,U_{i_{p(1)} i_{p(2)} i_{p(3)}} W_{j_1 j_2} \\
 & = &
\sum_{(j_1,j_2)\in{\mathfrak I}}\,\sum_{p\in\calP}\,C_{\epsilon , p}^{(j_1,j_2),\calP}(\nu)\,U_{i_{p(1)} i_{p(2)} i_{p(3)}} W_{j_1 j_2}\,
\end{eqnarray*}
where we use the notation
\begin{eqnarray} \label{equ5.16}%
C_{\epsilon , p}^{(j_1,j_2),\calP}(\nu) & := &
A_{\epsilon , p}^{(j_1,j_2)} + \sum_{s\in\calS_3\setminus\calP}\,A_{\epsilon , s}^{(j_1,j_2)}\,B_{s p}^{\calP}(\nu)\,.
\end{eqnarray}
If all rational functions $C_{\epsilon , p}^{(j_1,j_2),\calP}(\nu)$ fulfill
$C_{\epsilon , p}^{(j_1,j_2),\calP}(\nu)\not\equiv 0$, then we can consider the generic case of all such $\nu$ for which all
$C_{\epsilon , p}^{(j_1,j_2),\calP}(\nu)$ does not vanish. In this case we have
\begin{eqnarray*}
{\rm length}\,{\mathfrak P}_{ijklr}^{\rm red} & = &
|{\mathfrak I}|\times |\calP|\;=\; 10\times 2\;=\;20\,.
\end{eqnarray*}
Obviously, such a situation exists if $\epsilon = -1$.

However, a short look at (\ref{equ5.16}) shows that
some of the $C_{\epsilon , p}^{(j_1,j_2),\calP}(\nu)$ could satisfy a relation $C_{\epsilon , p}^{(j_1,j_2),\calP}(\nu)\equiv 0$ if the $A_{\epsilon , p}^{(j_1,j_2)}$ and $B_{s p}^{\calP}(\nu)$ fulfill certain conditions. A trivial example yields the case $\epsilon = 1$. In this case we read from the expression ${\textstyle\frac{1}{24}}(y_{t'}^{\ast}(U\otimes W))_{i j k l r}$ in Appendix A that $A_{\epsilon , p}^{(j_1,j_2)} = 0$ if $(j_1,j_2) = (i,j)$ or $(j_1,j_2) = (k,l)$. This leads to $C_{\epsilon , p}^{(j_1,j_2),\calP}(\nu)\equiv 0$ for $(j_1,j_2) = (i,j)$ or $(j_1,j_2) = (k,l)$ and to the reduction of ${\rm length}\,{\mathfrak P}_{ijklr}^{\rm red}$ from 20 to 16 summands.

Assume now that there would be two sets $\calP_1, \calP_2\subset\calS_3$ of two permutations to which different numbers of $C_{\epsilon , p}^{(j_1,j_2),{\calP}_i}(\nu)$ with $C_{\epsilon , p}^{(j_1,j_2),{\calP}_i}(\nu)\not\equiv 0$ belong. Then we would observe a change of ${\rm length}\,{\mathfrak P}_{ijklr}^{\rm red}$ even in the generic case if we pass from $\calP_1$ to $\calP_2$.
\end{Rem}
\vspace{5pt}

\subsubsection{The non-generic cases} \label{sec5.2.2}%
Now we handel non-generic cases of the search for short formulas for ${\rm length}\,{\mathfrak P}_{ijklr}^{\rm red}$. We start with the expressions (\ref{equ5.11}) which we found in the investigation of the generic case in Section \ref{sec5.2.1} and carry out the following
\begin{Proc} \label{proc5.4}%
Consider the formula (\ref{equ5.11}) which was won for a fixed set $\calP = \{p_1,p_2\}\subset\calS_3$ in the generic case.
\begin{enumerate}
\item{Determine the set $N_{\calP}$ of all $\nu\in\bbK$ for which at least one of the polynomials $P_q^{\calP}(\nu)$ from (\ref{equ5.11}) vanishes.}
\item{Remove all such $\nu\in N_{\calP}$ from $N_{\calP}$ that are roots of ${\Delta}_{\calP}(\nu)$. Further remove $\nu = \frac{1}{2}$ from $N_{\calP}$. (Compare Remark \ref{rem5.1a}.)}
\item{Set every $\nu\in N_{\calP}$ into the formula (\ref{equ5.11}) belonging to $\calP$ and determine the number of summands in the resulting expression.}
\end{enumerate}
\end{Proc}
The results of {\sf PERMS}-calculations according to Procedure \ref{proc5.4} are listed in Tables 3 and 4.
The {\sf Mathematica} notebooks of these calculations are \verb|roots12a.nb|, \ldots , \verb|roots56a.nb| in \cite{fie21}.
\begin{table}[t]
\begin{tabular}{c|c|c|c|c}
\hline
 & \multicolumn{2}{c|}{$\epsilon = 1$} & \multicolumn{2}{c}{$\epsilon = -1$}\\
$\calP$ & 
roots of $P_q^{\calP}(\nu)$ with & 
length of  & 
roots of $P_q^{\calP}(\nu)$ with  & 
length of \\
 & 
${\Delta}_{\calP}(\nu)\not= 0\;,\;\nu\not= 1/2$ & 
$\mathfrak{P}_{i_1\ldots i_5}^{\rm red}$ &  
${\Delta}_{\calP}(\nu)\not= 0\;,\;\nu\not= 1/2$ & 
$\mathfrak{P}_{i_1\ldots i_5}^{\rm red}$ \\
\hline
\hline
 12 &  0 & 12 &    & \\
    &  2 & 14 & 2  & 12 \\
\hline
 13 & -1 & 14 & -1 & 18 \\
\hline
 14 & -1 & 14 & -1 & 18 \\
    &  0 & 12 &    & \\
    &  2 & 12 &  2 & 10 \\
\hline
 15 & -1 & 12 & -1 & 10 \\
    &  0 & 12 &    & \\
    &  1 & 12 &    & \\
    &  2 & 14 &  2 & 12 \\
\hline
 16 & -1 & 12 & -1 & 10 \\
    &  0 & 12 &    & \\
    &  1 & 12 &    & \\
    &  2 & 12 &  2 & 10 \\
\hline
 23 & -1 & 14 & -1 & 18 \\
    &  0 & 12 &    & \\
    &  2 & 14 &  2 & 12 \\
\hline
 24 & -1 & 14 & -1 & 18 \\
    &  2 & 14 &  2 & 18 \\
\hline
 25 & -1 & 12 & -1 & 10 \\
    &  1 & 12 &    & \\
\hline
\hline
\end{tabular}
\vspace{3mm}
\caption{The lengths of $\mathfrak{P}_{i_1\ldots i_5}^{\rm red}$ for an $U$ from a ${\xi}_{\nu}$-symmetry class, where $\nu$ is an allowed root of a $P_q^{\calP}(\nu)$.}
\end{table}
\begin{table}[t]
\begin{tabular}{c|c|c|c|c}
\hline
 & \multicolumn{2}{c|}{$\epsilon = 1$} & \multicolumn{2}{c}{$\epsilon = -1$}\\
$\calP$ & 
roots of $P_q^{\calP}(\nu)$ with & 
length of  & 
roots of $P_q^{\calP}(\nu)$ with  & 
length of \\
 & 
${\Delta}_{\calP}(\nu)\not= 0\;,\;\nu\not= 1/2$ & 
$\mathfrak{P}_{i_1\ldots i_5}^{\rm red}$ &  
${\Delta}_{\calP}(\nu)\not= 0\;,\;\nu\not= 1/2$ & 
$\mathfrak{P}_{i_1\ldots i_5}^{\rm red}$ \\
\hline
\hline
 26 & -1 & 12 & -1 & 10 \\
    &  1 & 12 &    & \\
    &  2 & 14 &  2 & 18 \\
\hline
 34 &  0 & 12 &    & \\
    &  2 & 12 &  2 & 10 \\
\hline
 35 & -1 & 14 & -1 & 12 \\
    &  0 & 12 &    & \\
    &  1 & 12 &    & \\
    &  2 & 14 &  2 & 12 \\
\hline
 36 & -1 & 14 & -1 & 12 \\
    &  0 & 12 &    & \\
    &  1 & 12 &    & \\
    &  2 & 12 &  2 & 10 \\
\hline
 45 & -1 & 14 & -1 & 12 \\
    &  1 & 12 &    & \\
    &  2 & 14 &  2 & 18 \\
\hline
 46 & -1 & 14 & -1 & 12 \\
    &  1 & 12 & \\
\hline
 56 &  2 & 14 &  2 & 18 \\
\hline
\hline
\end{tabular}
\vspace{3mm}
\caption{Continuation: The lengths of $\mathfrak{P}_{i_1\ldots i_5}^{\rm red}$ for an $U$ from a ${\xi}_{\nu}$-symmetry class, where $\nu$ is an allowed root of a $P_q^{\calP}(\nu)$.}
\end{table}

A summary is given in Table (\ref{equ5.17}). This Table shows the minimal lengths of ${\mathfrak P}_{i_1\ldots i_5}^{\rm red}$ and the elements of the set
$N := \bigcup_{|\calP| = 2}\,N_{\calP}$ which is the collection of all (allowed) roots of the polynomials $P_q^{\calP}(\nu)$.
\begin{eqnarray} \label{equ5.17}%
 & &
\begin{array}{|c|c|c|}
\hline
 & \epsilon = 1 & \epsilon = -1 \\
\hline
\min {\rm length}\;{\mathfrak P}_{i_1\ldots i_5}^{\rm red} & 12 & 10 \\
\hline
N & \{ -1\,,\,0\,,\,1\,,\,2\} & \{ -1\,,\,2\} \\
\hline
\end{array}
\end{eqnarray}
For every $\nu\in N$ there is a set $\calP$ such that ${\mathfrak P}_{i_1\ldots i_5}^{\rm red}$ belonging to $\calP$ has the minimal length given in (\ref{equ5.17}). Note, however, that the roots $\nu = 0$ and $\nu = 1$ occur only for symmetric tensors $W$ (i.e. $\epsilon = 1$).

It is remarkable that for $\nu = -1$ and $\nu = 2$ sets $\calP$ can be found such that ${\mathfrak P}_{i_1\ldots i_5}^{\rm red}$ has minimal length both for
$\epsilon = 1$ and for $\epsilon = -1$. An example is the set $\calP = \{[1,2,3]\,,\,[3,2,1]\}$ denoted by ''16''. The expressions ${\mathfrak P}_{i_1\ldots i_5}^{\rm red}$ of minimal length belonging to ''16'' are presented in Appendix D. The example of the set ''16'' yields the proof of statement 3 in Theorem \ref{thm1.11}.\vspace{20pt}

\section{Index commutation symmetries} \label{sec6}%
In this section we investigate the question whether the tensors $U\in\calT_3 V$ belonging to $\eta$ or ${\xi}_{\nu}$ possess index commutation symmetries. We will see that such symmetries occur only for $\eta$ and a finite set of $\nu$-values. In particular the $\nu\in N$ given in (\ref{equ5.17}) are such $\nu$-values. This shows an interesting connection between index commutation symmetries and the maximal reduction of the length of ${\mathfrak P}_{i_1\ldots i_5}^{\rm red}$.
\begin{Def}
Let $C\subseteq\calS_r$ be a subgroup of $\calS_r$ and
$\theta : C \rightarrow \bbK^{\times}$ be a homomorphism of $C$ onto a finite subgroup of the multiplicative subgroup
$\bbK^{\times} := \bbK\setminus\{0\}$ of $\bbK$. We say that a tensor $T\in\calT_r V$ possesses the {\it index commutation symmetry} $(C,\theta)$ iff
\begin{eqnarray}
\forall\;c\in C\;:\;\;\;c\,T & = & \theta(c)\,T\,.
\end{eqnarray}
\end{Def}
We showed in \cite[p.115]{fie16}
\begin{Lem} \label{lemma6.2}%
A tensor $T\in\calT_r V$ possesses the symmetry $(C,\theta)$ iff $T$ satisfies
\begin{eqnarray}
{\tilde{\theta}}^{\ast}\,T & = & T
\end{eqnarray}
where
\begin{eqnarray} \label{equ6.3}%
\tilde{\theta} & := & \frac{1}{|C|}\,\sum_{c\in C}\,\theta(c)\,c
\end{eqnarray}
is the normalized symmetrizer belonging to $(C,\theta)$.
\end{Lem}
\begin{table}[t]
\begin{tabular}{c|c|l|c}
\hline
$C$ & $|C|$ & \parbox[c][16pt][c]{3cm}{\centerline{$|C|\cdot\tilde{\theta}$}} & name \\
\hline
 $\langle [2,1,3] \rangle$  & 2 & $[1,2,3] + [2,1,3]$  & sym2a \\
 $\langle [3,2,1] \rangle$  & 2 & $[1,2,3] + [3,2,1]$  & sym2b \\
 $\langle [1,3,2] \rangle$  & 2 & $[1,2,3] + [1,3,2]$  & sym2c \\
\hline
 $\langle [2,1,3] \rangle$  & 2 & $[1,2,3] - [2,1,3]$  & alt2a \\
 $\langle [3,2,1] \rangle$  & 2 & $[1,2,3] - [3,2,1]$  & alt2b \\
 $\langle [1,3,2] \rangle$  & 2 & $[1,2,3] - [1,3,2]$  & alt2c \\
\hline
 $\calA_3$  & 3 & $[1,2,3] + e^{2\imath\pi/3} [2,3,1] + e^{- 2\imath\pi/3} [3,1,2]$  & z1 \\
 $\calA_3$  & 3 & $[1,2,3] + e^{- 2\imath\pi/3} [2,3,1] + e^{2\imath\pi/3} [3,1,2]$  & z2 \\
 $\calA_3$  & 3 & $[1,2,3] + [2,3,1] + [3,1,2]$  & sym3 \\
\hline
 $\calS_3$  & 6 & $[1,2,3] + [1,3,2] + [2,1,3] +$  & sym6 \\
            &   & $[2,3,1] + [3,1,2] + [3,2,1]$  & \\
\hline
 $\calS_3$  & 6 & $[1,2,3] - [1,3,2] - [2,1,3] +$  & alt6 \\
            &   & $[2,3,1] + [3,1,2] - [3,2,1]$  & \\
\hline
\end{tabular}
\vspace{3mm}
\caption{All index commutation symmetries of subgroups of $\calS_3$}
\end{table}

Every symmetrizer (\ref{equ6.3}) is an idempotent (see \cite[p.115]{fie16}). Lemma \ref{lemma6.2} means that a tensor $T\in\calT_r V$ has a commutation symmetry $(C,\theta)$ iff $T$ belongs to the symmetry class defined by the right ideal $\tilde{\frakr} := {\tilde{\theta}}^{\ast}\cdot\bbK[\calS_r]$.
\begin{Def}
Let $\frakr\subseteq\bbK[\calS_r]$ be a right ideal defining a symmetry class of tensors $T\in\calT_r V$. We say that $\frakr$ admits a commutation symmetry $(C,\theta)$ over a subgroup $C\subseteq\calS_r$ iff
\begin{eqnarray} \label{equ6.4}%
\frakr & \subseteq & \tilde{\frakr} \;=\; {\tilde{\theta}}^{\ast}\cdot\bbK[\calS_r]\,.
\end{eqnarray}
\end{Def}
The validity of a condition (\ref{equ6.4}) can be checked comfortably by
\begin{Lem}\footnote{See H. Boerner \cite[pp.54-59]{boerner}.}
Let $e , \tilde{e}\in\bbK[\calS_r]$ be generating idempotents of right ideals
$\frakr , \tilde{\frakr}\subseteq\bbK[\calS_r]$. Then it holds
\begin{eqnarray}
\tilde{\frakr}\;\supseteq\;\frakr & \Leftrightarrow & \tilde{e}\cdot e\;=\;e\,.
\end{eqnarray}
\end{Lem}
Consequently, we have to check the conditions
\begin{eqnarray}
{\tilde{\theta}}^{\ast}\cdot\eta\;=\;\eta & {\rm or} &
{\tilde{\theta}}^{\ast}\cdot{\xi}_{\nu}\;=\;{\xi}_{\nu}
\end{eqnarray}
if we want to find out whether a tensor $U\in\calT_3 V$ from a symmetry class defined by $\eta$ or ${\xi}_{\nu}$ admits a commutation symmetry $(C,\theta)$, $C\subseteq\calS_3$.

In Table 5 we give a complete list of all index commutation symmetries that can be defined on arbitrary subgroups $C\subseteq\calS_3$. This list can be taken from the list \cite[pp.179]{fie16} of all index commutation symmetries belonging to subgroups of $\calS_r$ with $r\le 6$. (See also the {\sf HTML} version 
\begin{verbatim}
   home.t-online.de/home/Bernd.Fiedler.RoschStr.Leipzig/tensym.htm
\end{verbatim}
of the Table from \cite{fie16}.)

Now a computer calculation by means of {\sf PERMS} \cite{fie10} leads to the following Theorem \ref{thm6.5}. The {\sf Mathematica} notebook of the calculation can be found in \cite[comsym.nb]{fie21}.
\begin{Thm} \label{thm6.5}%
We consider the idempotents $\tilde{\theta}$ from Table {\rm 5} and the idempotents $e$ from the set $\{\eta\}\cup\{{\xi}_{\nu}\;|\;\nu\in\bbC\}$. Then a condition ${\tilde{\theta}}^{\ast}\cdot e = e$ is satisfied iff the pair $(\tilde{\theta}\,,\,e)$ is one of the pairs from the following table
\begin{center}
{\rm
\begin{tabular}{|c|l|}
\hline
\parbox[c][16pt][c]{1cm}{\centerline{$\tilde{\theta}$}} & $e$ \\
\hline
sym2a & $e = {\xi}_{\nu}$ with $\nu = 0$ \\
sym2b & $e = \eta$ \\
sym2c & $e = {\xi}_{\nu}$ with $\nu = 1$ \\
alt2a & $e = {\xi}_{\nu}$ with $\nu = 2$ \\
alt2b & $e = {\xi}_{\nu}$ with $\nu = \frac{1}{2}$ \\
alt2c & $e = {\xi}_{\nu}$ with $\nu = -1$ \\
z1 & $e = {\xi}_{\nu}$ with $\nu = e^{\imath\pi/3}$ \\
z2 & $e = {\xi}_{\nu}$ with $\nu = e^{- \imath\pi/3}$ \\
\hline
\end{tabular}
}
\end{center}
\end{Thm}\vspace{5pt}

We see that no of the considered tensors $U$ admits a commutation symmetry sym3, sym6 or alt6.
Furthermore the mentioned above connection between commutation symmetries and the maximal reduction of the length of ${\mathfrak P}_{i_1\ldots i_5}^{\rm red}$ becomes visible. If $\dim V\ge 3$, then we obtain statement 4 of Theorem \ref{thm1.11} from Theorem \ref{thm6.5} and Tables 3 and 4.\vspace{20pt}

\section{Tensors $U$ generated by covariant derivatives} \label{sec7}%
In \cite{fie03a} we showed that examples of tensors which lie in an irreducible symmetry class belonging to $\lambda = (2\,1)\vdash 3$ can be constructed from covariant derivatives of symmetric or alternating tensor fields of order 2. Such tensors are special examples of tensors $U$ considered in the present paper.
More precisely we proved the following facts in \cite{fie03a}:

Let $M$ be a differentiable\footnote{Here the word ''differentiable'' denotes the class $C^{\infty}$.} {\it manifold} of dimension $\dim M \ge 1$ and $\nabla$ be a {\it torsion-free covariant derivative} on $M$. Further let 
$\psi\in\calT_2 M$, $\omega\in\calT_2 M$ be covariant, differentiable tensor fields of order 2 which are {\it symmetric} or {\it skew-symmetric}, respectively. Then the tensors
\begin{eqnarray} \label{equ7.1}%
(\nabla\psi - {\rm sym}\,(\nabla\psi))|_p & \;\;\;,\;\;\; &
(\nabla\omega - {\rm alt}\,(\nabla\omega))|_p\;\;\;,\;\;\;p\in M
\end{eqnarray}
lie in irreducible symmetry classes\footnote{The symmetry classes of these tensor fields are different.} belonging to $\lambda = (2\,1)\vdash 3$. The operators 'sym' and 'alt' denote the {\it symmetrization} and {\it anti-symmetrization}, respectively. The right ideals that define the irreducible $(2\,1)$-symmetry classes for the tensors (\ref{equ7.1}) are generated by the idempotents
\begin{eqnarray}
h_s\;:=\;e_s - f_s & \;\;\;,\;\;\; &
h_a\;:=\;e_a - f_a\,
\end{eqnarray}
where
\begin{eqnarray}
e_s\;:=\;\frac{1}{2}\,([1,2,3] + [2,1,3]) & \;\;\;,\;\;\; &
e_a\;:=\;\frac{1}{2}\,([1,2,3] - [2,1,3]) \\
f_s\;:=\;\frac{1}{6}\,\sum_{p\in\calS_3}\,p & \;\;\;,\;\;\; &
f_a\;:=\;\frac{1}{6}\,\sum_{p\in\calS_3}\,{\rm sign}(p)\,p\,.
\end{eqnarray}
Note that this statement is based on the convention
\begin{eqnarray}
(\nabla T)_{i_1 i_2 i_3} & = & T_{i_1 i_2\,;\,i_3}\;=\;{\nabla}_{i_3}T_{i_1 i_2}
\end{eqnarray}
for the numbering of the tensor indices.

A second result of \cite{fie03a} was that the minimal right ideals
$h_s\cdot\bbK[\calS_3]$ and $h_a\cdot\bbK[\calS_3]$ are different from the critical right ideals $\frakr_0 = f\cdot\bbK[\calS_3]$ from Theorem \ref{thm1.10}. Thus the tensors (\ref{equ7.1}) can be used in generator formulas for algebraic covariant derivative curvature tensors.

Now we want to clear the question which of the idempotents from
\[
\{\eta\}\cup\{{\xi}_{\nu}\;|\;\nu\in\bbK\}
\]
generate the right ideals
$h_s\cdot\bbK[\calS_3]$ and $h_a\cdot\bbK[\calS_3]$. A short computer calculation \cite[derivs.nb]{fie21} by means of {\sf PERMS} \cite{fie10} yield
\begin{Thm}
The idempotents $h_s$, $h_a$, $\eta$, ${\xi}_{\nu}$ satisfy the following relations:
\begin{eqnarray}
\eta\cdot h_s\;\not=\;h_s & \;\;\;,\;\;\; & \eta\cdot h_a\;\not=\;h_a \\
{\xi}_{\nu}\cdot h_s\;=\;h_s & \Leftrightarrow & \nu\;=\;0 \\
{\xi}_{\nu}\cdot h_a\;=\;h_a & \Leftrightarrow & \nu\;=\;2\,.
\end{eqnarray}
Thus we obtain
\begin{eqnarray}
{\xi}_0\cdot\bbK[\calS_3]\;=\;h_s\cdot\bbK[\calS_3] & \;\;\;,\;\;\; &
{\xi}_2\cdot\bbK[\calS_3]\;=\;h_a\cdot\bbK[\calS_3]\,.
\end{eqnarray}
\end{Thm}

It is interesting that $\nu = 0$ and $\nu = 2$ are $\nu$-values which allow the construction of shortest formulas for ${\mathfrak P}_{i_1\ldots i_5}^{\rm red}$ (see Tables 3 and 4). Note, however, that only $\nu = 2$ lead to tensors $U$ which can be used to construct such shortest formulas both for symmetric and for skew-symmetric $W$. A tensor $U$ with $\nu = 0$ will not produce a minimal length of ${\mathfrak P}_{i_1\ldots i_5}^{\rm red}$ if $W$ is an alternating tensor. If we now assume that $\dim M\ge 3$, then the equation system (\ref{equ4.4}) yields a complete set of linear identities for the tensors (\ref{equ7.1}) (see explanation of (\ref{equ4.1a})) and we obtain Theorem \ref{thm1.12}.
\clearpage

\section*{Appendix A: First reduction of $\mathfrak{P}_{i_1\ldots i_5}$}
In Appendix A we present tensor coordinates of the tensors (\ref{equ1.22}) which we calculated by means of the {\sf Mathematica} packages {\sf PERMS} \cite{fie10} and {\sf Ricci} \cite{ricci3}. {\sf PERMS} contains a tool \verb|Symmetrize| which makes it possible for us to apply symmetry operators $a\in\bbK[\calS_r]$ of {\sf Perms} onto the coordinates $T_{i_1\ldots i_r}$ of tensors $T\in\calT_r V$ defined in {\sf Ricci} and to form
\begin{eqnarray*}
(aT)_{i_1\ldots i_r} & = & \sum_{p\in\calS_r}\,a(p)\,T_{i_{p(1)}\ldots i_{p(r)}}\,.
\end{eqnarray*}

We use the common symbol $W_{i j}$ for the coordinates $A_{i j}$ and $S_{i j}$ of the tensors $A, S\in\calT_2 V$. Then we apply the idempotent $\frac{1}{24}\,y_{t'}^{\ast}\in\bbK[\calS_5]$ to the coordinates $U_{i j k} W_{l r}$ or $W_{i j} U_{k l r}$ of the tensors $U\otimes W$ or $W\otimes U$, respectively. The resulting expressions
\begin{eqnarray*}
{\textstyle\frac{1}{24}}(y_{t'}^{\ast}(U\otimes W))_{i j k l r} & = &
{\textstyle\frac{1}{24}}\left\{U_{i j k} W_{l r} - U_{i j l} W_{k r} -
U_{i j r} W_{k l} + U_{i j r} W_{l k} \pm 44\,\mathrm{terms} \right\}\\
{\textstyle\frac{1}{24}}(y_{t'}^{\ast}(W\otimes U))_{i j k l r} & = &
{\textstyle\frac{1}{24}}\left\{W_{i j} U_{k l r} - W_{i j} U_{l k r} -
W_{i j} U_{r k l} + W_{i j} U_{r l k} \pm 44\,\mathrm{terms} \right\}
\end{eqnarray*}
are automatically reduced by {\sf Ricci} by means of the identity
\begin{eqnarray*}
W_{i j} & = & \epsilon\,W_{j i}\;\;\;,\;\;\;\epsilon\in\{ 1\,,\, -1 \}\,.
\end{eqnarray*}
We could get a further reduction of these coordinate expressions if we considered all linear identities which are fulfilled by the tensors of the symmetry class of $U$. However, we want to carry out such a reduction only in the following appendices. During the calculations of Appendix A we assume that $U$ is a tensor ''without any symmetry'', i.e. $U\in\calT_3 V$ is a tensor whose symmetry class is defined by the right ideal $\frakr = \bbK[\calS_3]$.

Under these assumptions we obtain the following two coordinate expressions for
$\frac{1}{24}(y_{t'}^{\ast}(U\otimes W))_{i j k l r}$ and $\frac{1}{24}(y_{t'}^{\ast}(W\otimes U))_{i j k l r}$ (see \cite[part1.nb]{fie21}):
\clearpage

\subsubsection*{The coordinates $\frac{1}{24}(y_{t'}^{\ast}(U\otimes W))_{i j k l r}$}
\begin{eqnarray*}
 & &
\frac{-1 + \epsilon}{24} \,\tensor{U}{\down{k}\down{l}\down{r}} 
    \tensor{W}{\down{i}\down{j}} + 
  \frac{1 - \epsilon}{24} \,\tensor{U}{\down{l}\down{k}\down{r}} 
    \tensor{W}{\down{i}\down{j}} + 
  \frac{1 - \epsilon}{24} \,\tensor{U}{\down{r}\down{k}\down{l}} 
    \tensor{W}{\down{i}\down{j}} + 
  \frac{-1 + \epsilon}{24} \,\tensor{U}{\down{r}\down{l}\down{k}} 
    \tensor{W}{\down{i}\down{j}} - \\ & &
  \frac{1}{24} \,\tensor{U}{\down{j}\down{l}\down{r}} 
    \tensor{W}{\down{i}\down{k}} - 
  \frac{\epsilon}{24} \,\tensor{U}{\down{l}\down{j}\down{r}} 
    \tensor{W}{\down{i}\down{k}} - 
  \frac{\epsilon}{24} \,\tensor{U}{\down{r}\down{j}\down{l}} 
    \tensor{W}{\down{i}\down{k}} - 
  \frac{1}{24} \,\tensor{U}{\down{r}\down{l}\down{j}} 
    \tensor{W}{\down{i}\down{k}} + \\ & &
  \frac{1}{24} \,\tensor{U}{\down{j}\down{k}\down{r}} 
    \tensor{W}{\down{i}\down{l}} + 
  \frac{\epsilon}{24} \,\tensor{U}{\down{k}\down{j}\down{r}} 
    \tensor{W}{\down{i}\down{l}} + 
  \frac{\epsilon}{24} \,\tensor{U}{\down{r}\down{j}\down{k}} 
    \tensor{W}{\down{i}\down{l}} + 
  \frac{1}{24} \,\tensor{U}{\down{r}\down{k}\down{j}} 
    \tensor{W}{\down{i}\down{l}} + \\ & &
  \frac{1}{24} \,\tensor{U}{\down{j}\down{k}\down{l}} 
    \tensor{W}{\down{i}\down{r}} - 
  \frac{1}{24} \,\tensor{U}{\down{j}\down{l}\down{k}} 
    \tensor{W}{\down{i}\down{r}} - 
  \frac{1}{24} \,\tensor{U}{\down{k}\down{l}\down{j}} 
    \tensor{W}{\down{i}\down{r}} + 
  \frac{1}{24} \,\tensor{U}{\down{l}\down{k}\down{j}} 
    \tensor{W}{\down{i}\down{r}} + \\ & &
  \frac{1}{24} \,\tensor{U}{\down{i}\down{l}\down{r}} 
    \tensor{W}{\down{j}\down{k}} + 
  \frac{\epsilon}{24} \,\tensor{U}{\down{l}\down{i}\down{r}} 
    \tensor{W}{\down{j}\down{k}} + 
  \frac{\epsilon}{24} \,\tensor{U}{\down{r}\down{i}\down{l}} 
    \tensor{W}{\down{j}\down{k}} + 
  \frac{1}{24} \,\tensor{U}{\down{r}\down{l}\down{i}} 
    \tensor{W}{\down{j}\down{k}} - \\ & &
  \frac{1}{24} \,\tensor{U}{\down{i}\down{k}\down{r}} 
    \tensor{W}{\down{j}\down{l}} - 
  \frac{\epsilon}{24} \,\tensor{U}{\down{k}\down{i}\down{r}} 
    \tensor{W}{\down{j}\down{l}} - 
  \frac{\epsilon}{24} \,\tensor{U}{\down{r}\down{i}\down{k}} 
    \tensor{W}{\down{j}\down{l}} - 
  \frac{1}{24} \,\tensor{U}{\down{r}\down{k}\down{i}} 
    \tensor{W}{\down{j}\down{l}} - \\ & &
  \frac{1}{24} \,\tensor{U}{\down{i}\down{k}\down{l}} 
    \tensor{W}{\down{j}\down{r}} + 
  \frac{1}{24} \,\tensor{U}{\down{i}\down{l}\down{k}} 
    \tensor{W}{\down{j}\down{r}} + 
  \frac{1}{24} \,\tensor{U}{\down{k}\down{l}\down{i}} 
    \tensor{W}{\down{j}\down{r}} - 
  \frac{1}{24} \,\tensor{U}{\down{l}\down{k}\down{i}} 
    \tensor{W}{\down{j}\down{r}} + \\ & &
  \frac{-1 + \epsilon}{24} \,\tensor{U}{\down{i}\down{j}\down{r}} 
    \tensor{W}{\down{k}\down{l}} + 
  \frac{1 - \epsilon}{24} \,\tensor{U}{\down{j}\down{i}\down{r}} 
    \tensor{W}{\down{k}\down{l}} + 
  \frac{1 - \epsilon}{24} \,\tensor{U}{\down{r}\down{i}\down{j}} 
    \tensor{W}{\down{k}\down{l}} + 
  \frac{-1 + \epsilon}{24} \,\tensor{U}{\down{r}\down{j}\down{i}} 
    \tensor{W}{\down{k}\down{l}} - \\ & &
  \frac{1}{24} \,\tensor{U}{\down{i}\down{j}\down{l}} 
    \tensor{W}{\down{k}\down{r}} + 
  \frac{1}{24} \,\tensor{U}{\down{j}\down{i}\down{l}} 
    \tensor{W}{\down{k}\down{r}} + 
  \frac{1}{24} \,\tensor{U}{\down{l}\down{i}\down{j}} 
    \tensor{W}{\down{k}\down{r}} - 
  \frac{1}{24} \,\tensor{U}{\down{l}\down{j}\down{i}} 
    \tensor{W}{\down{k}\down{r}} + \\ & &
  \frac{1}{24} \,\tensor{U}{\down{i}\down{j}\down{k}} 
    \tensor{W}{\down{l}\down{r}} - 
  \frac{1}{24} \,\tensor{U}{\down{j}\down{i}\down{k}} 
    \tensor{W}{\down{l}\down{r}} - 
  \frac{1}{24} \,\tensor{U}{\down{k}\down{i}\down{j}} 
    \tensor{W}{\down{l}\down{r}} + 
  \frac{1}{24} \,\tensor{U}{\down{k}\down{j}\down{i}} 
    \tensor{W}{\down{l}\down{r}}
\end{eqnarray*}

\subsubsection*{The coordinates $\frac{1}{24}(y_{t'}^{\ast}(W\otimes U))_{i j k l r}$}
\begin{eqnarray*}
 & &
\frac{1 - \epsilon}{24} \,\tensor{U}{\down{k}\down{l}\down{r}} 
    \tensor{W}{\down{i}\down{j}} + 
  \frac{-1 + \epsilon}{24} \,\tensor{U}{\down{l}\down{k}\down{r}} 
    \tensor{W}{\down{i}\down{j}} + 
  \frac{-1 + \epsilon}{24} \,\tensor{U}{\down{r}\down{k}\down{l}} 
    \tensor{W}{\down{i}\down{j}} + 
  \frac{1 - \epsilon}{24} \,\tensor{U}{\down{r}\down{l}\down{k}} 
    \tensor{W}{\down{i}\down{j}} - \\ & &
  \frac{\epsilon}{24} \,\tensor{U}{\down{j}\down{l}\down{r}} 
    \tensor{W}{\down{i}\down{k}} - 
  \frac{1}{24} \,\tensor{U}{\down{l}\down{j}\down{r}} 
    \tensor{W}{\down{i}\down{k}} - 
  \frac{1}{24} \,\tensor{U}{\down{r}\down{j}\down{l}} 
    \tensor{W}{\down{i}\down{k}} - 
  \frac{\epsilon}{24} \,\tensor{U}{\down{r}\down{l}\down{j}} 
    \tensor{W}{\down{i}\down{k}} + \\ & &
  \frac{\epsilon}{24} \,\tensor{U}{\down{j}\down{k}\down{r}} 
    \tensor{W}{\down{i}\down{l}} + 
  \frac{1}{24} \,\tensor{U}{\down{k}\down{j}\down{r}} 
    \tensor{W}{\down{i}\down{l}} + 
  \frac{1}{24} \,\tensor{U}{\down{r}\down{j}\down{k}} 
    \tensor{W}{\down{i}\down{l}} + 
  \frac{\epsilon}{24} \,\tensor{U}{\down{r}\down{k}\down{j}} 
    \tensor{W}{\down{i}\down{l}} + \\ & &
  \frac{\epsilon}{24} \,\tensor{U}{\down{j}\down{k}\down{l}} 
    \tensor{W}{\down{i}\down{r}} - 
  \frac{\epsilon}{24} \,\tensor{U}{\down{j}\down{l}\down{k}} 
    \tensor{W}{\down{i}\down{r}} - 
  \frac{\epsilon}{24} \,\tensor{U}{\down{k}\down{l}\down{j}} 
    \tensor{W}{\down{i}\down{r}} + 
  \frac{\epsilon}{24} \,\tensor{U}{\down{l}\down{k}\down{j}} 
    \tensor{W}{\down{i}\down{r}} + \\ & &
  \frac{\epsilon}{24} \,\tensor{U}{\down{i}\down{l}\down{r}} 
    \tensor{W}{\down{j}\down{k}} + 
  \frac{1}{24} \,\tensor{U}{\down{l}\down{i}\down{r}} 
    \tensor{W}{\down{j}\down{k}} + 
  \frac{1}{24} \,\tensor{U}{\down{r}\down{i}\down{l}} 
    \tensor{W}{\down{j}\down{k}} + 
  \frac{\epsilon}{24} \,\tensor{U}{\down{r}\down{l}\down{i}} 
    \tensor{W}{\down{j}\down{k}} - \\ & &
  \frac{\epsilon}{24} \,\tensor{U}{\down{i}\down{k}\down{r}} 
    \tensor{W}{\down{j}\down{l}} - 
  \frac{1}{24} \,\tensor{U}{\down{k}\down{i}\down{r}} 
    \tensor{W}{\down{j}\down{l}} - 
  \frac{1}{24} \,\tensor{U}{\down{r}\down{i}\down{k}} 
    \tensor{W}{\down{j}\down{l}} - 
  \frac{\epsilon}{24} \,\tensor{U}{\down{r}\down{k}\down{i}} 
    \tensor{W}{\down{j}\down{l}} - \\ & & 
  \frac{\epsilon}{24} \,\tensor{U}{\down{i}\down{k}\down{l}} 
    \tensor{W}{\down{j}\down{r}} + 
  \frac{\epsilon}{24} \,\tensor{U}{\down{i}\down{l}\down{k}} 
    \tensor{W}{\down{j}\down{r}} + 
  \frac{\epsilon}{24} \,\tensor{U}{\down{k}\down{l}\down{i}} 
    \tensor{W}{\down{j}\down{r}} - 
  \frac{\epsilon}{24} \,\tensor{U}{\down{l}\down{k}\down{i}} 
    \tensor{W}{\down{j}\down{r}} + \\ & &
  \frac{1 - \epsilon}{24} \,\tensor{U}{\down{i}\down{j}\down{r}} 
    \tensor{W}{\down{k}\down{l}} + 
  \frac{-1 + \epsilon}{24} \,\tensor{U}{\down{j}\down{i}\down{r}} 
    \tensor{W}{\down{k}\down{l}} + 
  \frac{-1 + \epsilon}{24} \,\tensor{U}{\down{r}\down{i}\down{j}} 
    \tensor{W}{\down{k}\down{l}} + 
  \frac{1 - \epsilon}{24} \,\tensor{U}{\down{r}\down{j}\down{i}} 
    \tensor{W}{\down{k}\down{l}} - \\ & & 
  \frac{\epsilon}{24} \,\tensor{U}{\down{i}\down{j}\down{l}} 
    \tensor{W}{\down{k}\down{r}} + 
  \frac{\epsilon}{24} \,\tensor{U}{\down{j}\down{i}\down{l}} 
    \tensor{W}{\down{k}\down{r}} + 
  \frac{\epsilon}{24} \,\tensor{U}{\down{l}\down{i}\down{j}} 
    \tensor{W}{\down{k}\down{r}} - 
  \frac{\epsilon}{24} \,\tensor{U}{\down{l}\down{j}\down{i}} 
    \tensor{W}{\down{k}\down{r}} + \\ & &
  \frac{\epsilon}{24} \,\tensor{U}{\down{i}\down{j}\down{k}} 
    \tensor{W}{\down{l}\down{r}} - 
  \frac{\epsilon}{24} \,\tensor{U}{\down{j}\down{i}\down{k}} 
    \tensor{W}{\down{l}\down{r}} - 
  \frac{\epsilon}{24} \,\tensor{U}{\down{k}\down{i}\down{j}} 
    \tensor{W}{\down{l}\down{r}} + 
  \frac{\epsilon}{24} \,\tensor{U}{\down{k}\down{j}\down{i}} 
    \tensor{W}{\down{l}\down{r}}
\end{eqnarray*}
$\,$\\*[0.3cm]
\noindent Furthermore our calculations in the notebook \cite[part1.nb]{fie21} yield the following\\*[0.3cm]
{\bf Lemma:} {\it The above tensor coordinates satisfy}
\begin{eqnarray*}
{\textstyle\frac{1}{24}}(y_{t'}^{\ast}(U\otimes W))_{i j k l r} & = &
\epsilon {\textstyle\frac{1}{24}}(y_{t'}^{\ast}(W\otimes U))_{i j k l r}
\end{eqnarray*}
{\it if} $W_{ji} = \epsilon W_{ij}$, $\epsilon\in\{1,-1\}$.  {\it The expressions possesses {\rm 32} summands if $W$ is symmetric and {\rm 40} summands if $W$ is an alternating tensor.}\vspace{20pt}

\section*{Appendix B: Shortest $\mathfrak{P}_{i_1\ldots i_5}^{\rm red}$ with $U$ belonging to $\eta$}
In Appendix B we present formulas of minimal length for $\frac{1}{24}(y_{t'}^{\ast}(U\otimes W))_{i j k l r}$ where $U$ is a tensor from the symmetry class defined by the idempotent $\eta$ from (\ref{equ3.21}). Such formulas of minimal lenght arise if $\calP = \{[1,2,3]\,,\,[1,3,2]\}$. Linear identities which characterize the symmetry of $U$ in this case are given in (\ref{equ4.11}).

We obtain a length of 12 summands if $W$ is symmetric ($\epsilon = 1$) and a length of 20 summands if $W$ is an alternating tensor ($\epsilon = -1$). The determination of the formulas is described in Section \ref{sec5.1}. The computer calculations for Appendix B can be found in the {\sf Mathematica} notebook \cite[part12b.nb]{fie21}.

\subsubsection*{The coordinates $\frac{1}{24}(y_{t'}^{\ast}(U\otimes W))_{i j k l r}$ for $\calP = \{[1,2,3]\,,\,[1,3,2]\}$, $\epsilon = 1$}
\begin{eqnarray*}
 & &
\frac{1}{12} \,\tensor{U}{\down{j}\down{r}\down{l}} 
    \tensor{W}{\down{i}\down{k}} - 
  \frac{1}{12} \,\tensor{U}{\down{j}\down{r}\down{k}} 
    \tensor{W}{\down{i}\down{l}} + 
  \frac{1}{12} \,\tensor{U}{\down{j}\down{k}\down{l}} 
    \tensor{W}{\down{i}\down{r}} - 
  \frac{1}{12} \,\tensor{U}{\down{j}\down{l}\down{k}} 
    \tensor{W}{\down{i}\down{r}} - \\ & &
  \frac{1}{12} \,\tensor{U}{\down{i}\down{r}\down{l}} 
    \tensor{W}{\down{j}\down{k}} + 
  \frac{1}{12} \,\tensor{U}{\down{i}\down{r}\down{k}} 
    \tensor{W}{\down{j}\down{l}} - 
  \frac{1}{12} \,\tensor{U}{\down{i}\down{k}\down{l}} 
    \tensor{W}{\down{j}\down{r}} + 
  \frac{1}{12} \,\tensor{U}{\down{i}\down{l}\down{k}} 
    \tensor{W}{\down{j}\down{r}} - \\ & &
  \frac{1}{6} \,\tensor{U}{\down{i}\down{j}\down{l}} 
    \tensor{W}{\down{k}\down{r}} - 
  \frac{1}{12} \,\tensor{U}{\down{i}\down{l}\down{j}} 
    \tensor{W}{\down{k}\down{r}} + 
  \frac{1}{6} \,\tensor{U}{\down{i}\down{j}\down{k}} 
    \tensor{W}{\down{l}\down{r}} + 
  \frac{1}{12} \,\tensor{U}{\down{i}\down{k}\down{j}} 
    \tensor{W}{\down{l}\down{r}}
\end{eqnarray*}

\subsubsection*{The coordinates $\frac{1}{24}(y_{t'}^{\ast}(U\otimes W))_{i j k l r}$ for $\calP = \{[1,2,3]\,,\,[1,3,2]\}$, $\epsilon = -1$}
\begin{eqnarray*}
 & &
- \frac{1}{3} \,\tensor{U}{\down{k}\down{l}\down{r}} 
    \tensor{W}{\down{i}\down{j}} - 
  \frac{1}{6} \,\tensor{U}{\down{k}\down{r}\down{l}} 
    \tensor{W}{\down{i}\down{j}} - 
  \frac{1}{6} \,\tensor{U}{\down{j}\down{l}\down{r}} 
    \tensor{W}{\down{i}\down{k}} - 
  \frac{1}{12} \,\tensor{U}{\down{j}\down{r}\down{l}} 
    \tensor{W}{\down{i}\down{k}} + \\ & &
  \frac{1}{6} \,\tensor{U}{\down{j}\down{k}\down{r}} 
    \tensor{W}{\down{i}\down{l}} + 
  \frac{1}{12} \,\tensor{U}{\down{j}\down{r}\down{k}} 
    \tensor{W}{\down{i}\down{l}} + 
  \frac{1}{12} \,\tensor{U}{\down{j}\down{k}\down{l}} 
    \tensor{W}{\down{i}\down{r}} - 
  \frac{1}{12} \,\tensor{U}{\down{j}\down{l}\down{k}} 
    \tensor{W}{\down{i}\down{r}} + \\ & &
  \frac{1}{6} \,\tensor{U}{\down{i}\down{l}\down{r}} 
    \tensor{W}{\down{j}\down{k}} + 
  \frac{1}{12} \,\tensor{U}{\down{i}\down{r}\down{l}} 
    \tensor{W}{\down{j}\down{k}} - 
  \frac{1}{6} \,\tensor{U}{\down{i}\down{k}\down{r}} 
    \tensor{W}{\down{j}\down{l}} - 
  \frac{1}{12} \,\tensor{U}{\down{i}\down{r}\down{k}} 
    \tensor{W}{\down{j}\down{l}} - \\ & &
  \frac{1}{12} \,\tensor{U}{\down{i}\down{k}\down{l}} 
    \tensor{W}{\down{j}\down{r}} + 
  \frac{1}{12} \,\tensor{U}{\down{i}\down{l}\down{k}} 
    \tensor{W}{\down{j}\down{r}} - 
  \frac{1}{3} \,\tensor{U}{\down{i}\down{j}\down{r}} 
    \tensor{W}{\down{k}\down{l}} - 
  \frac{1}{6} \,\tensor{U}{\down{i}\down{r}\down{j}} 
    \tensor{W}{\down{k}\down{l}} - \\ & &
  \frac{1}{6} \,\tensor{U}{\down{i}\down{j}\down{l}} 
    \tensor{W}{\down{k}\down{r}} - 
  \frac{1}{12} \,\tensor{U}{\down{i}\down{l}\down{j}} 
    \tensor{W}{\down{k}\down{r}} + 
  \frac{1}{6} \,\tensor{U}{\down{i}\down{j}\down{k}} 
    \tensor{W}{\down{l}\down{r}} + 
  \frac{1}{12} \,\tensor{U}{\down{i}\down{k}\down{j}} 
    \tensor{W}{\down{l}\down{r}}
\end{eqnarray*}
\vspace{20pt}

\section*{Appendix C: Generic case for $\mathfrak{P}_{i_1\ldots i_5}^{\rm red}$ with $U$ belonging to ${\xi}_{\nu}$}
In Appendix C we present formulas of minimal length for $\frac{1}{24}(y_{t'}^{\ast}(U\otimes W))_{i j k l r}$ where $U$ is a tensor from the symmetry class defined by the idempotent ${\xi}_{\nu}$ from (\ref{equ3.20}).

First we consider the example of the set $\calP = \{[1,2,3]\,,\,[1,3,2]\}$ and restrict us to such $\nu\in\bbK$ for which ${\Delta}_{\calP}(\nu)\not= 0$. Linear identities which characterize the symmetry of $U$ in this case are given in (\ref{equ5.13}).

We obtain a length of 16 summands if $W$ is symmetric ($\epsilon = 1$) and a length of 20 summands if $W$ is an alternating tensor ($\epsilon = -1$). The determination of the formulas is described in Section \ref{sec5.2}. The computer calculations can be found in the {\sf Mathematica} notebook \cite[part12a.nb]{fie21}.

\subsection*{$\mathfrak{P}_{i_1\ldots i_5}^{\rm red}$ for $\calP = \{[1,2,3]\,,\,[1,3,2]\}$}

\subsubsection*{The coordinates $\frac{1}{24}(y_{t'}^{\ast}(U\otimes W))_{i j k l r}$ for $\epsilon = 1$}
\begin{eqnarray*}
 & &
\frac{-\left( -1 + 2\,\nu \right) }{24\,\left( -1 + \nu \right) \,\left( 1 + \nu \right) } 
   \,\,\tensor{U}{\down{j}\down{l}\down{r}} 
    \tensor{W}{\down{i}\down{k}} + 
  \frac{\nu\,\left( -1 + 2\,\nu \right) }
    {24\,\left( -1 + \nu \right) \,\left( 1 + \nu \right) } 
   \,\,\tensor{U}{\down{j}\down{r}\down{l}} 
    \tensor{W}{\down{i}\down{k}} + \\ & &
  \frac{-1 + 2\,\nu}{24\,\left( -1 + \nu \right) \,\left( 1 + \nu \right) } 
   \,\,\tensor{U}{\down{j}\down{k}\down{r}} 
    \tensor{W}{\down{i}\down{l}} - 
  \frac{\nu\,\left( -1 + 2\,\nu \right) }
    {24\,\left( -1 + \nu \right) \,\left( 1 + \nu \right) } 
   \,\,\tensor{U}{\down{j}\down{r}\down{k}} 
    \tensor{W}{\down{i}\down{l}} + \\ & &
  \frac{-1 + 2\,\nu}{24\,\left( -1 + \nu \right) } 
   \,\,\tensor{U}{\down{j}\down{k}\down{l}} 
    \tensor{W}{\down{i}\down{r}} - 
  \frac{-1 + 2\,\nu}{24\,\left( -1 + \nu \right) } 
   \,\,\tensor{U}{\down{j}\down{l}\down{k}} 
    \tensor{W}{\down{i}\down{r}} + \\ & &
  \frac{-1 + 2\,\nu}{24\,\left( -1 + \nu \right) \,\left( 1 + \nu \right) } 
   \,\,\tensor{U}{\down{i}\down{l}\down{r}} 
    \tensor{W}{\down{j}\down{k}} - 
  \frac{\nu\,\left( -1 + 2\,\nu \right) }
    {24\,\left( -1 + \nu \right) \,\left( 1 + \nu \right) } 
   \,\,\tensor{U}{\down{i}\down{r}\down{l}} 
    \tensor{W}{\down{j}\down{k}} - \\ & &
  \frac{-1 + 2\,\nu}{24\,\left( -1 + \nu \right) \,\left( 1 + \nu \right) } 
   \,\,\tensor{U}{\down{i}\down{k}\down{r}} 
    \tensor{W}{\down{j}\down{l}} + 
  \frac{\nu\,\left( -1 + 2\,\nu \right) }
    {24\,\left( -1 + \nu \right) \,\left( 1 + \nu \right) } 
   \,\,\tensor{U}{\down{i}\down{r}\down{k}} 
    \tensor{W}{\down{j}\down{l}} - \\ & &
  \frac{-1 + 2\,\nu}{24\,\left( -1 + \nu \right) } 
   \,\,\tensor{U}{\down{i}\down{k}\down{l}} 
    \tensor{W}{\down{j}\down{r}} + 
  \frac{-1 + 2\,\nu}{24\,\left( -1 + \nu \right) } 
   \,\,\tensor{U}{\down{i}\down{l}\down{k}} 
    \tensor{W}{\down{j}\down{r}} - \\ & &
  \frac{{\left( -1 + 2\,\nu \right) }^2}
    {24\,\left( -1 + \nu \right) \,\left( 1 + \nu \right) } 
   \,\,\tensor{U}{\down{i}\down{j}\down{l}} 
    \tensor{W}{\down{k}\down{r}} - 
  \frac{\left( -2 + \nu \right) \,\left( -1 + 2\,\nu \right) }
    {24\,\left( -1 + \nu \right) \,\left( 1 + \nu \right) } 
   \,\,\tensor{U}{\down{i}\down{l}\down{j}} 
    \tensor{W}{\down{k}\down{r}} + \\ & &
  \frac{{\left( -1 + 2\,\nu \right) }^2}
    {24\,\left( -1 + \nu \right) \,\left( 1 + \nu \right) } 
   \,\,\tensor{U}{\down{i}\down{j}\down{k}} 
    \tensor{W}{\down{l}\down{r}} + 
  \frac{\left( -2 + \nu \right) \,\left( -1 + 2\,\nu \right) }
    {24\,\left( -1 + \nu \right) \,\left( 1 + \nu \right) } 
   \,\,\tensor{U}{\down{i}\down{k}\down{j}} 
    \tensor{W}{\down{l}\down{r}}
\end{eqnarray*}
\clearpage

\subsubsection*{The coordinates $\frac{1}{24}(y_{t'}^{\ast}(U\otimes W))_{i j k l r}$ for $\epsilon = -1$}
\begin{eqnarray*}
 & &
\frac{-{\left( -1 + 2\,\nu \right) }^2}
    {12\,\left( -1 + \nu \right) \,\left( 1 + \nu \right) } 
   \,\,\tensor{U}{\down{k}\down{l}\down{r}} 
    \tensor{W}{\down{i}\down{j}} - 
  \frac{\left( -2 + \nu \right) \,\left( -1 + 2\,\nu \right) }
    {12\,\left( -1 + \nu \right) \,\left( 1 + \nu \right) } 
   \,\,\tensor{U}{\down{k}\down{r}\down{l}} 
    \tensor{W}{\down{i}\down{j}} - \\ & &
  \frac{{\left( -1 + 2\,\nu \right) }^2}
    {24\,\left( -1 + \nu \right) \,\left( 1 + \nu \right) } 
   \,\,\tensor{U}{\down{j}\down{l}\down{r}} 
    \tensor{W}{\down{i}\down{k}} - 
  \frac{\left( -2 + \nu \right) \,\left( -1 + 2\,\nu \right) }
    {24\,\left( -1 + \nu \right) \,\left( 1 + \nu \right) } 
   \,\,\tensor{U}{\down{j}\down{r}\down{l}} 
    \tensor{W}{\down{i}\down{k}} + \\ & &
  \frac{{\left( -1 + 2\,\nu \right) }^2}
    {24\,\left( -1 + \nu \right) \,\left( 1 + \nu \right) } 
   \,\,\tensor{U}{\down{j}\down{k}\down{r}} 
    \tensor{W}{\down{i}\down{l}} + 
  \frac{\left( -2 + \nu \right) \,\left( -1 + 2\,\nu \right) }
    {24\,\left( -1 + \nu \right) \,\left( 1 + \nu \right) } 
   \,\,\tensor{U}{\down{j}\down{r}\down{k}} 
    \tensor{W}{\down{i}\down{l}} + \\ & &
  \frac{-1 + 2\,\nu}{24\,\left( -1 + \nu \right) } 
   \,\,\tensor{U}{\down{j}\down{k}\down{l}} 
    \tensor{W}{\down{i}\down{r}} - 
  \frac{-1 + 2\,\nu}{24\,\left( -1 + \nu \right) } 
   \,\,\tensor{U}{\down{j}\down{l}\down{k}} 
    \tensor{W}{\down{i}\down{r}} + \\ & &
  \frac{{\left( -1 + 2\,\nu \right) }^2}
    {24\,\left( -1 + \nu \right) \,\left( 1 + \nu \right) } 
   \,\,\tensor{U}{\down{i}\down{l}\down{r}} 
    \tensor{W}{\down{j}\down{k}} + 
  \frac{\left( -2 + \nu \right) \,\left( -1 + 2\,\nu \right) }
    {24\,\left( -1 + \nu \right) \,\left( 1 + \nu \right) } 
   \,\,\tensor{U}{\down{i}\down{r}\down{l}} 
    \tensor{W}{\down{j}\down{k}} - \\ & &
  \frac{{\left( -1 + 2\,\nu \right) }^2}
    {24\,\left( -1 + \nu \right) \,\left( 1 + \nu \right) } 
   \,\,\tensor{U}{\down{i}\down{k}\down{r}} 
    \tensor{W}{\down{j}\down{l}} - 
  \frac{\left( -2 + \nu \right) \,\left( -1 + 2\,\nu \right) }
    {24\,\left( -1 + \nu \right) \,\left( 1 + \nu \right) } 
   \,\,\tensor{U}{\down{i}\down{r}\down{k}} 
    \tensor{W}{\down{j}\down{l}} - \\ & &
  \frac{-1 + 2\,\nu}{24\,\left( -1 + \nu \right) } 
   \,\,\tensor{U}{\down{i}\down{k}\down{l}} 
    \tensor{W}{\down{j}\down{r}} + 
  \frac{-1 + 2\,\nu}{24\,\left( -1 + \nu \right) } 
   \,\,\tensor{U}{\down{i}\down{l}\down{k}} 
    \tensor{W}{\down{j}\down{r}} - \\ & &
  \frac{{\left( -1 + 2\,\nu \right) }^2}
    {12\,\left( -1 + \nu \right) \,\left( 1 + \nu \right) } 
   \,\,\tensor{U}{\down{i}\down{j}\down{r}} 
    \tensor{W}{\down{k}\down{l}} - 
  \frac{\left( -2 + \nu \right) \,\left( -1 + 2\,\nu \right) }
    {12\,\left( -1 + \nu \right) \,\left( 1 + \nu \right) } 
   \,\,\tensor{U}{\down{i}\down{r}\down{j}} 
    \tensor{W}{\down{k}\down{l}} - \\ & &
  \frac{{\left( -1 + 2\,\nu \right) }^2}
    {24\,\left( -1 + \nu \right) \,\left( 1 + \nu \right) } 
   \,\,\tensor{U}{\down{i}\down{j}\down{l}} 
    \tensor{W}{\down{k}\down{r}} - 
  \frac{\left( -2 + \nu \right) \,\left( -1 + 2\,\nu \right) }
    {24\,\left( -1 + \nu \right) \,\left( 1 + \nu \right) } 
   \,\,\tensor{U}{\down{i}\down{l}\down{j}} 
    \tensor{W}{\down{k}\down{r}} + \\ & &
  \frac{{\left( -1 + 2\,\nu \right) }^2}
    {24\,\left( -1 + \nu \right) \,\left( 1 + \nu \right) } 
   \,\,\tensor{U}{\down{i}\down{j}\down{k}} 
    \tensor{W}{\down{l}\down{r}} + 
  \frac{\left( -2 + \nu \right) \,\left( -1 + 2\,\nu \right) }
    {24\,\left( -1 + \nu \right) \,\left( 1 + \nu \right) } 
   \,\,\tensor{U}{\down{i}\down{k}\down{j}} 
    \tensor{W}{\down{l}\down{r}}
\end{eqnarray*}
$\,$\\*[0.2cm]
{\bf Remark:} We see that no denominator of the coefficients in the above expressions has the root $\nu = \frac{1}{2}$ both for $\epsilon = 1$ and for $\epsilon = -1$. Furthermore every numerator of the above coefficients contains a factor $(-1 + 2\nu)$. Consequently we obtain $\frac{1}{24}(y_{t'}^{\ast}(U\otimes W))_{i j k l r}|_{\nu = 1/2} = 0$ both for $\epsilon = 1$ and for $\epsilon = -1$. This illustrates that a tensor $U$ from the symmetry class of $\xi_{1/2}$ can not be used to generate a non-trivial algebraic covariant derivative curvature tensor.\vspace{0.5cm}

Now we repeat the above considerations for the set $\calP = \{[1,2,3]\,,\,[3,2,1]\}$. For this set $\calP$ the polynomial ${\Delta}_{\calP}(\nu)$ possesses the only root $\nu = \frac{1}{2}$, which is the critical $\nu$-value for our construction of algebraic covariant derivative curvature tensors from the tensors $U$ and $W$.

Since no other $\nu$-values has to be excluded in this case, the formulas for $\mathfrak{P}_{i_1\ldots i_5}^{\rm red}$ are formulas which yield algebraic covariant derivative curvature tensors for every $\nu\not=\frac{1}{2}$ (see Remark \ref{rem5.1}). The computer calculations can be found in the {\sf Mathematica} notebook \cite[part16a.nb]{fie21}.

For $\calP = \{[1,2,3]\,,\,[3,2,1]\}$ Procedure \ref{proc4.6} yields the following linear identities (\ref{equ4.10a}) for $U$:
\begin{eqnarray*}
 & &
\begin{array}{ccccccccccc}
- & \frac{{\nu}^2 - \nu + 1}{2{\nu} - 1}\,U_{ijk} & + & \frac{{\nu}^2 - 1}{2{\nu} - 1}\,U_{kji}        & +  &         &         &         & U_{ikj} & = & 0 \\
  & \frac{{\nu}^2 - \nu + 1}{2{\nu} - 1}\,U_{ijk} & - & \frac{{\nu}^2 - 2\nu}{2{\nu} - 1}\,U_{kji} & + &         &         & U_{jik} &         & = & 0 \\
-  & \frac{{\nu}^2 - 2\nu}{2{\nu} - 1}\,U_{ijk}        & + & \frac{{\nu}^2 - \nu + 1}{2{\nu} - 1}\,U_{kji} & + &         & U_{jki} &         &         & = & 0 \\
  & \frac{{\nu}^2 - 1}{2{\nu} - 1}\,U_{ijk} & - & \frac{{\nu}^2 - \nu + 1}{2{\nu} - 1}\,U_{kji} & + & U_{kij} &         &         &         & = & 0 \\
\end{array}\;.
\end{eqnarray*}
If we carry out the Procedure \ref{proc4.8} by means of these identities, we obtain the following expressions for $\mathfrak{P}_{i_1\ldots i_5}^{\rm red}$:

\subsection*{$\mathfrak{P}_{i_1\ldots i_5}^{\rm red}$ for $\calP = \{[1,2,3]\,,\,[3,2,1]\}$}

\subsubsection*{The coordinates $\frac{1}{24}(y_{t'}^{\ast}(U\otimes W))_{i j k l r}$ for $\epsilon = 1$}
\begin{eqnarray*}
 & &
\frac{-1 + \nu}{24} \,\tensor{U}{\down{j}\down{l}\down{r}} 
    \tensor{W}{\down{i}\down{k}} - 
  \frac{\nu}{24} \,\tensor{U}{\down{r}\down{l}\down{j}} 
    \tensor{W}{\down{i}\down{k}} + 
  \frac{1 - \nu}{24} \,\tensor{U}{\down{j}\down{k}\down{r}} 
    \tensor{W}{\down{i}\down{l}} + 
  \frac{\nu}{24} \,\tensor{U}{\down{r}\down{k}\down{j}} 
    \tensor{W}{\down{i}\down{l}} + \\ & &
  \frac{2 - \nu}{24} \,\tensor{U}{\down{j}\down{k}\down{l}} 
    \tensor{W}{\down{i}\down{r}} + 
  \frac{1 + \nu}{24} \,\tensor{U}{\down{l}\down{k}\down{j}} 
    \tensor{W}{\down{i}\down{r}} + 
  \frac{1 - \nu}{24} \,\tensor{U}{\down{i}\down{l}\down{r}} 
    \tensor{W}{\down{j}\down{k}} + 
  \frac{\nu}{24} \,\tensor{U}{\down{r}\down{l}\down{i}} 
    \tensor{W}{\down{j}\down{k}} + \\ & &
  \frac{-1 + \nu}{24} \,\tensor{U}{\down{i}\down{k}\down{r}} 
    \tensor{W}{\down{j}\down{l}} - 
  \frac{\nu}{24} \,\tensor{U}{\down{r}\down{k}\down{i}} 
    \tensor{W}{\down{j}\down{l}} + 
  \frac{-2 + \nu}{24} \,\tensor{U}{\down{i}\down{k}\down{l}} 
    \tensor{W}{\down{j}\down{r}} + 
  \frac{-1 - \nu}{24} \,\tensor{U}{\down{l}\down{k}\down{i}} 
    \tensor{W}{\down{j}\down{r}} + \\ & &
  \frac{-1 - \nu}{24} \,\tensor{U}{\down{i}\down{j}\down{l}} 
    \tensor{W}{\down{k}\down{r}} + 
  \frac{-2 + \nu}{24} \,\tensor{U}{\down{l}\down{j}\down{i}} 
    \tensor{W}{\down{k}\down{r}} + 
  \frac{1 + \nu}{24} \,\tensor{U}{\down{i}\down{j}\down{k}} 
    \tensor{W}{\down{l}\down{r}} + 
  \frac{2 - \nu}{24} \,\tensor{U}{\down{k}\down{j}\down{i}} 
    \tensor{W}{\down{l}\down{r}}
\end{eqnarray*}

\subsubsection*{The coordinates $\frac{1}{24}(y_{t'}^{\ast}(U\otimes W))_{i j k l r}$ for $\epsilon = -1$}
\begin{eqnarray*}
 & &
\frac{-1 - \nu}{12} \,\tensor{U}{\down{k}\down{l}\down{r}} 
    \tensor{W}{\down{i}\down{j}} + 
  \frac{-2 + \nu}{12} \,\tensor{U}{\down{r}\down{l}\down{k}} 
    \tensor{W}{\down{i}\down{j}} + 
  \frac{-1 - \nu}{24} \,\tensor{U}{\down{j}\down{l}\down{r}} 
    \tensor{W}{\down{i}\down{k}} + 
  \frac{-2 + \nu}{24} \,\tensor{U}{\down{r}\down{l}\down{j}} 
    \tensor{W}{\down{i}\down{k}} + \\ & &
  \frac{1 + \nu}{24} \,\tensor{U}{\down{j}\down{k}\down{r}} 
    \tensor{W}{\down{i}\down{l}} + 
  \frac{2 - \nu}{24} \,\tensor{U}{\down{r}\down{k}\down{j}} 
    \tensor{W}{\down{i}\down{l}} + 
  \frac{2 - \nu}{24} \,\tensor{U}{\down{j}\down{k}\down{l}} 
    \tensor{W}{\down{i}\down{r}} + 
  \frac{1 + \nu}{24} \,\tensor{U}{\down{l}\down{k}\down{j}} 
    \tensor{W}{\down{i}\down{r}} + \\ & &
  \frac{1 + \nu}{24} \,\tensor{U}{\down{i}\down{l}\down{r}} 
    \tensor{W}{\down{j}\down{k}} + 
  \frac{2 - \nu}{24} \,\tensor{U}{\down{r}\down{l}\down{i}} 
    \tensor{W}{\down{j}\down{k}} + 
  \frac{-1 - \nu}{24} \,\tensor{U}{\down{i}\down{k}\down{r}} 
    \tensor{W}{\down{j}\down{l}} + 
  \frac{-2 + \nu}{24} \,\tensor{U}{\down{r}\down{k}\down{i}} 
    \tensor{W}{\down{j}\down{l}} + \\ & &
  \frac{-2 + \nu}{24} \,\tensor{U}{\down{i}\down{k}\down{l}} 
    \tensor{W}{\down{j}\down{r}} + 
  \frac{-1 - \nu}{24} \,\tensor{U}{\down{l}\down{k}\down{i}} 
    \tensor{W}{\down{j}\down{r}} + 
  \frac{-1 - \nu}{12} \,\tensor{U}{\down{i}\down{j}\down{r}} 
    \tensor{W}{\down{k}\down{l}} + 
  \frac{-2 + \nu}{12} \,\tensor{U}{\down{r}\down{j}\down{i}} 
    \tensor{W}{\down{k}\down{l}} + \\ & &
  \frac{-1 - \nu}{24} \,\tensor{U}{\down{i}\down{j}\down{l}} 
    \tensor{W}{\down{k}\down{r}} + 
  \frac{-2 + \nu}{24} \,\tensor{U}{\down{l}\down{j}\down{i}} 
    \tensor{W}{\down{k}\down{r}} + 
  \frac{1 + \nu}{24} \,\tensor{U}{\down{i}\down{j}\down{k}} 
    \tensor{W}{\down{l}\down{r}} + 
  \frac{2 - \nu}{24} \,\tensor{U}{\down{k}\down{j}\down{i}} 
    \tensor{W}{\down{l}\down{r}}
\end{eqnarray*}
$\,$\\*[0.2cm]
{\bf Remark:} During the above calculation the critical factor $(-1 + 2\nu)$ was canceled in all fractions within the expressions for $\frac{1}{24}(y_{t'}^{\ast}(U\otimes W))_{i j k l r}$. Thus we could set $\nu = \frac{1}{2}$ in these formulas. However it is clear that the resulting formulas do not describe an algebraic covariant derivative curvature tensor, because the assumption $\nu\not=\frac{1}{2}$ was the foundation of our calculation. The above formulas represent an algebraic covariant derivative curvature tensor only then if the coordinates of $U$ also satisfy the above linear identities of type (\ref{equ4.10a}). However, these identities are not defined for $\nu = \frac{1}{2}$.\vspace{20pt}

\section*{Appendix D: Non-generic cases for $\mathfrak{P}_{i_1\ldots i_5}^{\rm red}$ with $U$ belonging to ${\xi}_{\nu}$}
Now we present examples for a further reduction of the length of the formulas for $\frac{1}{24}(y_{t'}^{\ast}(U\otimes W))_{i j k l r}$ from Appendix C. Again we assume that $U$ is a tensor from the symmetry class defined by the idempotent ${\xi}_{\nu}$ from (\ref{equ3.20}). However we use such a value $\nu$ which is a root of a polynomial $P_q^{\calP}(\nu)$ occuring in the considered expression from Appendix C. The meaning of $P_q^{\calP}(\nu)$ is defined in (\ref{equ5.11}). The vanishing of $P_q^{\calP}(\nu)$ leads to a reduction of the length of the considered expression of type (\ref{equ5.11}).

Our consideration is based on the set $\calP = \{[1,2,3]\,,\,[3,2,1]\}$. For this set the values $\nu = -1$ and $\nu = 2$ lead to the minimal length of  $\mathfrak{P}_{i_1\ldots i_5}^{\rm red}$ both for symmetric $W$ ($\epsilon = 1$) and for alternating $W$ ($\epsilon = 1$). We obtain a length of 12 summands if $W$ is symmetric ($\epsilon = 1$) and a length of 10 summands if $W$ is an alternating tensor ($\epsilon = -1$). The determination of these results is described in Section \ref{sec5.2.2}. The computer calculations can be found in the {\sf Mathematica} notebook \cite[roots16a.nb]{fie21}.

The below expressions are coordinates of algebraic covariant derivative curvature tensors only if the coordinates of $U$ satisfy the correct linear identities describing the symmetry of $U$. These identities can be obtained if one sets $\nu = -1$ or $\nu = 2$ into the general identities for $U$ given in the second part of Appendix C.

\subsection*{$\mathfrak{P}_{i_1\ldots i_5}^{\rm red}$ containing a $U$ from the symmetry class defined by ${\xi}_{-1}$}

\subsubsection*{The coordinates $\frac{1}{24}(y_{t'}^{\ast}(U\otimes W))_{i j k l r}$ for $\epsilon = 1$}
\begin{eqnarray*}
 & &
-\frac{1}{12}  \,\tensor{U}{\down{j}\down{l}\down{r}} 
    \tensor{W}{\down{i}\down{k}} + 
  \frac{1}{24} \,\tensor{U}{\down{r}\down{l}\down{j}} 
    \tensor{W}{\down{i}\down{k}} + 
  \frac{1}{12} \,\tensor{U}{\down{j}\down{k}\down{r}} 
    \tensor{W}{\down{i}\down{l}} - 
  \frac{1}{24} \,\tensor{U}{\down{r}\down{k}\down{j}} 
    \tensor{W}{\down{i}\down{l}} + 
  \frac{1}{8} \,\tensor{U}{\down{j}\down{k}\down{l}} 
    \tensor{W}{\down{i}\down{r}} + \\ & &
  \frac{1}{12} \,\tensor{U}{\down{i}\down{l}\down{r}} 
    \tensor{W}{\down{j}\down{k}} - 
  \frac{1}{24} \,\tensor{U}{\down{r}\down{l}\down{i}} 
    \tensor{W}{\down{j}\down{k}} - 
  \frac{1}{12} \,\tensor{U}{\down{i}\down{k}\down{r}} 
    \tensor{W}{\down{j}\down{l}} + 
  \frac{1}{24} \,\tensor{U}{\down{r}\down{k}\down{i}} 
    \tensor{W}{\down{j}\down{l}} - 
  \frac{1}{8} \,\tensor{U}{\down{i}\down{k}\down{l}} 
    \tensor{W}{\down{j}\down{r}} - \\ & &
  \frac{1}{8} \,\tensor{U}{\down{l}\down{j}\down{i}} 
    \tensor{W}{\down{k}\down{r}} + 
  \frac{1}{8} \,\tensor{U}{\down{k}\down{j}\down{i}} 
    \tensor{W}{\down{l}\down{r}}
\end{eqnarray*}

\subsubsection*{The coordinates $\frac{1}{24}(y_{t'}^{\ast}(U\otimes W))_{i j k l r}$ for $\epsilon = -1$}
\begin{eqnarray*}
 & &
-\frac{1}{4}  \,\tensor{U}{\down{r}\down{l}\down{k}} 
    \tensor{W}{\down{i}\down{j}} - 
  \frac{1}{8} \,\tensor{U}{\down{r}\down{l}\down{j}} 
    \tensor{W}{\down{i}\down{k}} + 
  \frac{1}{8} \,\tensor{U}{\down{r}\down{k}\down{j}} 
    \tensor{W}{\down{i}\down{l}} + 
  \frac{1}{8} \,\tensor{U}{\down{j}\down{k}\down{l}} 
    \tensor{W}{\down{i}\down{r}} + 
  \frac{1}{8} \,\tensor{U}{\down{r}\down{l}\down{i}} 
    \tensor{W}{\down{j}\down{k}} - \\ & &
  \frac{1}{8} \,\tensor{U}{\down{r}\down{k}\down{i}} 
    \tensor{W}{\down{j}\down{l}} - 
  \frac{1}{8} \,\tensor{U}{\down{i}\down{k}\down{l}} 
    \tensor{W}{\down{j}\down{r}} - 
  \frac{1}{4} \,\tensor{U}{\down{r}\down{j}\down{i}} 
    \tensor{W}{\down{k}\down{l}} - 
  \frac{1}{8} \,\tensor{U}{\down{l}\down{j}\down{i}} 
    \tensor{W}{\down{k}\down{r}} + 
  \frac{1}{8} \,\tensor{U}{\down{k}\down{j}\down{i}} 
    \tensor{W}{\down{l}\down{r}}
\end{eqnarray*}

\subsection*{$\mathfrak{P}_{i_1\ldots i_5}^{\rm red}$ containing a $U$ from the symmetry class defined by ${\xi}_2$}

\subsubsection*{The coordinates $\frac{1}{24}(y_{t'}^{\ast}(U\otimes W))_{i j k l r}$ for $\epsilon = 1$}
\begin{eqnarray*}
 & &
\frac{1}{24} \,\tensor{U}{\down{j}\down{l}\down{r}} 
    \tensor{W}{\down{i}\down{k}} - 
  \frac{1}{12} \,\tensor{U}{\down{r}\down{l}\down{j}} 
    \tensor{W}{\down{i}\down{k}} - 
  \frac{1}{24} \,\tensor{U}{\down{j}\down{k}\down{r}} 
    \tensor{W}{\down{i}\down{l}} + 
  \frac{1}{12} \,\tensor{U}{\down{r}\down{k}\down{j}} 
    \tensor{W}{\down{i}\down{l}} + 
  \frac{1}{8} \,\tensor{U}{\down{l}\down{k}\down{j}} 
    \tensor{W}{\down{i}\down{r}} - \\ & &
  \frac{1}{24} \,\tensor{U}{\down{i}\down{l}\down{r}} 
    \tensor{W}{\down{j}\down{k}} + 
  \frac{1}{12} \,\tensor{U}{\down{r}\down{l}\down{i}} 
    \tensor{W}{\down{j}\down{k}} + 
  \frac{1}{24} \,\tensor{U}{\down{i}\down{k}\down{r}} 
    \tensor{W}{\down{j}\down{l}} - 
  \frac{1}{12} \,\tensor{U}{\down{r}\down{k}\down{i}} 
    \tensor{W}{\down{j}\down{l}} - 
  \frac{1}{8} \,\tensor{U}{\down{l}\down{k}\down{i}} 
    \tensor{W}{\down{j}\down{r}} - \\ & &
  \frac{1}{8} \,\tensor{U}{\down{i}\down{j}\down{l}} 
    \tensor{W}{\down{k}\down{r}} + 
  \frac{1}{8} \,\tensor{U}{\down{i}\down{j}\down{k}} 
    \tensor{W}{\down{l}\down{r}}
\end{eqnarray*}

\subsubsection*{The coordinates $\frac{1}{24}(y_{t'}^{\ast}(U\otimes W))_{i j k l r}$ for $\epsilon = -1$}
\begin{eqnarray*}
 & &
-\frac{1}{4}  \,\tensor{U}{\down{k}\down{l}\down{r}} 
    \tensor{W}{\down{i}\down{j}} - 
  \frac{1}{8} \,\tensor{U}{\down{j}\down{l}\down{r}} 
    \tensor{W}{\down{i}\down{k}} + 
  \frac{1}{8} \,\tensor{U}{\down{j}\down{k}\down{r}} 
    \tensor{W}{\down{i}\down{l}} + 
  \frac{1}{8} \,\tensor{U}{\down{l}\down{k}\down{j}} 
    \tensor{W}{\down{i}\down{r}} + 
  \frac{1}{8} \,\tensor{U}{\down{i}\down{l}\down{r}} 
    \tensor{W}{\down{j}\down{k}} - \\ & &
  \frac{1}{8} \,\tensor{U}{\down{i}\down{k}\down{r}} 
    \tensor{W}{\down{j}\down{l}} - 
  \frac{1}{8} \,\tensor{U}{\down{l}\down{k}\down{i}} 
    \tensor{W}{\down{j}\down{r}} - 
  \frac{1}{4} \,\tensor{U}{\down{i}\down{j}\down{r}} 
    \tensor{W}{\down{k}\down{l}} - 
  \frac{1}{8} \,\tensor{U}{\down{i}\down{j}\down{l}} 
    \tensor{W}{\down{k}\down{r}} + 
  \frac{1}{8} \,\tensor{U}{\down{i}\down{j}\down{k}} 
    \tensor{W}{\down{l}\down{r}}
\end{eqnarray*}\vspace{0.5cm}

\noindent {\bf Acknowledgements.} I would like to thank Prof. P. B. Gilkey for important and helpful discussions and for valuable suggestions for future investigations.
\vspace{0.5cm}

\begin{thebibliography}{10}

\bibitem{bbg97}
N.~Bla{\v z}i{\'c}, N.~Bokan and P.~Gilkey.
\newblock A note on {O}sserman {L}orentzian manifolds.
\newblock {\em Bull. London Math. Soc.}, 29: 227--230, 1997.

\bibitem{boerner}
H.~Boerner.
\newblock {\em Darstellungen von Gruppen}, volume~74 of {\em Die Grundlehren
  der mathematischen Wissenschaften in Einzeldarstellungen}.
\newblock Springer-Verlag, Berlin, G\"{o}ttingen, Heidelberg, 1955.

\bibitem{boerner2}
H.~Boerner.
\newblock {\em Representations of Groups}.
\newblock North-Holland Publishing Company, Amsterdam, 2. revised edition,
  1970.

\bibitem{c88}
Q.-S. Chi.
\newblock A curvature characterization of certain locally rank-one symmetric
  spaces.
\newblock {\em J. Differential Geom.}, 28: 187--202, 1988.

\bibitem{clausbaum1}
M.~Clausen and U.~Baum.
\newblock {\em Fast Fourier Transforms}.
\newblock BI Wissenschaftsverlag, Mannheim, Leipzig, Wien, Z\"urich, 1993.

\bibitem{clausbaum2}
M.~Clausen and U.~Baum.
\newblock Fast {F}ourier transforms for symmetric groups.
\newblock In L.~Finkelstein and W.~M. Kantor, editors, {\em Groups and
  Computation}, volume~11 of {\em DIMACS Series in Discrete Mathematics and
  Theoretical Computer Science}, pages 27 -- 39, Providence, Rhode Island,
  1993. DIMACS, American Mathematical Society.
\newblock Proc. DIMACS workshop Rutgers University 1991.

\bibitem{fie21}
B.~Fiedler.
\newblock Examples of calculations by means of {PERMS}. {M}athematica
  notebooks.
\newblock Internet
  \verb|http://home.t-online.de/home/Bernd.Fiedler.RoschStr.Leipzig/pnbks.htm|.

\bibitem{fie12}
B.~Fiedler.
\newblock A characterization of the dependence of the {R}iemannian metric on
  the curvature tensor by {Y}oung symmetrizers.
\newblock {\em Z. Anal. Anw.}, 17(1): 135 -- 157, 1998.

\bibitem{fie16}
B.~Fiedler.
\newblock {\em An Algorithm for the Decomposition of Ideals of Semi-Simple
  Rings and its Application to Symbolic Tensor Calculations by Computer}.
\newblock Habilitationsschrift, Universit{\"a}t Leipzig, Leipzig, Germany,
  November 1999.
\newblock Fakult{\"a}t f{\"u}r Mathematik und Informatik.

\bibitem{fie10}
B.~Fiedler.
\newblock {\em {PERMS 2.1 (15.1.1999)}}.
\newblock Mathematisches Institut, Universit\"at Leipzig, Leipzig, 1999.
\newblock Will be sent in to MathSource, Wolfram Research Inc.

\bibitem{fie17}
B.~Fiedler.
\newblock Characterization of tensor symmetries by group ring subspaces and
  computation of normal forms of tensor coordinates.
\newblock In A.~Betten, A.~Kohnert, R.~Laue and A.~Wassermann,
  editors, {\em Algebraic Combinatorics and Applications. Proceedings of the
  Euroconference, ALCOMA, G{\"o}{\ss}weinstein, Germany, September 12--19,
  1999}, pages 118--133, Berlin, 2001. Springer-Verlag.

\bibitem{fie18}
B.~Fiedler.
\newblock Ideal decompositions and computation of tensor normal forms.
\newblock In {\em S\'eminaire Lotharingien de Combinatoire}, 2001.
\newblock El. published: \verb|http://www.mat.univie.ac.at/~slc|.
  B45g, 16 pp. Archive: \verb|http://arXiv.org/abs/math.CO/0211156|.

\bibitem{fie20}
B.~Fiedler.
\newblock Determination of the structure of algebraic curvature tensors by
  means of {Y}oung symmetrizers.
\newblock In {\em S\'eminaire Lotharingien de Combinatoire}, 2002.
\newblock Electronically published: \verb|http://www.mat.univie.ac.at/~slc|.
  B48d, 20 pp. Preprint: \verb|http://arXiv.org/abs/math.CO/0212278|.

\bibitem{fie03b}
B.~Fiedler.
\newblock Generators of algebraic covariant derivative curvature tensors and
  {Y}oung symmetrizers.
\newblock Preprint: \verb|http://arXiv.org/abs/math.CO/0310020|., 2003.
\newblock 18 pages. Chapter for a book "Progress in Computer Science Research",
  in preparation by Nova Science Publishers, Inc.

\bibitem{fie03a}
B.~Fiedler.
\newblock On the symmetry classes of the first covariant derivatives of tensor
  fields.
\newblock In {\em S\'eminaire Lotharingien de Combinatoire}, 2003.
\newblock Accepted by the electronic journal. 21 pp. Preprint:
  \verb|http://arXiv.org/abs/math.CO/0301042|.

\bibitem{full4}
S.~A. Fulling, R.~C. King, B.~G. Wybourne and C.~J. Cummins.
\newblock Normal forms for tensor polynomials: I. {T}he {R}iemann tensor.
\newblock {\em Class. Quantum Grav.}, 9: 1151 -- 1197, 1992.

\bibitem{fultharr}
W.~Fulton and J.~Harris.
\newblock {\em Representation Theory: A First Course}, volume 129 of {\em
  Graduate Texts in Mathematics, Readings in Mathematics}.
\newblock Springer-Verlag, New York, Berlin, Heidelberg, London, Paris, Tokyo,
  Hong Kong, Barcelona, Budapest, 1991.

\bibitem{fulton}
William Fulton.
\newblock {\em Young Tableaux}.
\newblock Number~35 in London Mathematical Society Student Texts. Cambridge
  University Press, Cambridge, New York, Melbourne, 1997.

\bibitem{gkv97}
E.~Garc{\'i}a-Ri{\'o}, D.~N. Kupeli and M.~E. V{\'a}zquez-Abal.
\newblock On a problem of {O}sserman in {L}orentzian geometry.
\newblock {\em Differential Geom. Appl.}, 7: 85--100, 1997.

\bibitem{gilkey3}
P.~B. Gilkey.
\newblock Geometric properties of the curvature operator.
\newblock In W.~H. Chen, A.-M. Li, U.~Simon, L.~Verstraelen, C.~P. Wang and
  M.~Wiehe, editors, {\em Proceedings of Beijing Conference ''Geometry and
  Topology of Submanifolds''}, volume~10, pages 62--70, Singapore, 2000. World
  Scientific.

\bibitem{gilkey5}
P.~B. Gilkey.
\newblock {\em Geometric Properties of Natural Operators Defined by the Riemann
  Curvature Tensor}.
\newblock World Scientific Publishing Co., Singapore, New Jersey, London, Hong
  Kong, 2001.

\bibitem{jameskerb}
G.~D. James and A.~Kerber.
\newblock {\em The Representation Theory of the Symmetric Group}, volume~16 of
  {\em Encyclopedia of Mathematics and its Applications}.
\newblock Addison-Wesley Publishing Company, Reading, Mass., London, Amsterdam,
  Don Mills, Ont., Sidney, Tokyo, 1981.

\bibitem{kerber}
A.~Kerber.
\newblock {\em Representations of Permutation Groups}, volume 240, 495 of {\em
  Lecture Notes in Mathematics}.
\newblock Springer-Verlag, Berlin, Heidelberg, New York, 1971, 1975.

\bibitem{kerber3}
A.~Kerber.
\newblock {\em Algebraic combinatorics via finite group actions}.
\newblock BI-Wiss.-Verl., Mannheim, Wien, Z\"urich, 1991.

\bibitem{kerbkohn2}
A.~Kerber and A.~Kohnert.
\newblock {\em {SYMMETRICA 1.0}}.
\newblock Lehrstuhl II f\"ur Mathematik, Department of Mathematics, University
  of Bayreuth, Bayreuth, 1994.
\newblock Manual.

\bibitem{kerbkohnlas}
A.~Kerber, A.~Kohnert and A.~Lascoux.
\newblock {SYMMETRICA}, an object oriented computer-algebra system for the
  symmetric group.
\newblock {\em J. Symbolic Computation}, 14: 195 -- 203, 1992.

\bibitem{ricci3}
J.~M. Lee, D.~Lear, J.~Roth, J.~Coskey and L.~Nave.
\newblock {\em Ricci. A Mathematica package for doing tensor calculations in
  differential geometry. User's Manual. Version 1.32}.
\newblock Department of Mathematics, Box 354350, University of Washington,
  Seattle, WA 98195-4350, 1992 - 1998.
\newblock Ricci's home page: \verb|http://www.math.washington.edu/~lee/Ricci/|.

\bibitem{littlew1}
D.~E. Littlewood.
\newblock {\em The Theory of Group Characters and Matrix Representations of
  Groups}.
\newblock Clarendon Press, Oxford, 2. edition, 1950.

\bibitem{mcdonald}
I.~G. Macdonald.
\newblock {\em Symmetric Functions and Hall Polynomials}.
\newblock Clarendon Press, Oxford, 1979.

\bibitem{muell}
W.~M\"{u}ller.
\newblock {\em Darstellungstheorie von endlichen Gruppen}.
\newblock Teubner Studienb\"{u}cher Mathematik. B. G. Teubner, Stuttgart, 1980.

\bibitem{naimark}
M.~A. Naimark and A.~I. \v{S}tern.
\newblock {\em Theory of Group Representations}, volume 246 of {\em Grundlehren
  der Mathematischen Wissenschaften}.
\newblock Springer-Verlag, Berlin, Heidelberg, New York, 1982.

\bibitem{n02}
Y.~Nikolayevsky.
\newblock Osserman conjecture in dimension $n \ne 8, 16$.
\newblock preprint, 2002.
\newblock http://arXiv.org/abs/math.DG/0204258.

\bibitem{oss90}
R.~Osserman.
\newblock Curvature in the eighties.
\newblock {\em Amer. Math. Monthly}, 97: 731--756, 1990.

\bibitem{waerden}
B.~L. {van der Waerden}.
\newblock {\em Algebra}, volume I, II.
\newblock Springer-Verlag, Berlin, Heidelberg, New York, London, Paris, Tokyo,
  Hong Kong, Barcelona, Budapest, 9., 6. edition, 1993.

\bibitem{weyl1}
H.~Weyl.
\newblock {\em The Classical Groups, their Invariants and Representations}.
\newblock Princeton University Press, Princeton, New Jersey, 1939.

\end{thebibliography}

\end{document}